\documentclass[12pt,amscd]{amsart}
\usepackage{amssymb}

\numberwithin{equation}{section}
\newtheorem{thm}[equation]{Theorem}

\newtheorem{lem}[equation]{Lemma}

\newtheorem{prop}[equation]{Proposition}

\newenvironment{pf}{\proof[\proofname]}{\endproof}
\newenvironment{pf*}[1]{\proof[#1]}{\endproof}

\newcommand{\comment}[1]{}

\theoremstyle{definition}

\newtheorem{defn}[equation]{Definition}

\theoremstyle{remark}

\newtheorem*{rmk}{Remark}
\newtheorem*{rmks}{Remarks}

\newtheorem*{ack}{Acknowledgement}

\makeatother

\begin{document}
\baselineskip=18truept


\def\C {{\mathbb C}}
\def\Cn {{\mathbb C}^n}
\def\R {{\mathbb R}}
\def\Rn {{\mathbb R}^n}
\def\Z {{\mathbb Z}}
\def\N {{\mathbb N}}
\def\cal#1{{\mathcal #1}}
\def\bb#1{{\mathbb #1}}

\def\dbar {\bar \partial }
\def\dir {{\mathcal D}}
\def\lev#1{{\mathcal L}{\left({#1}\right)}}
\def\lap {\Delta }

\def\ol {{\mathcal O}}
\def\ideal#1{{\mathcal I}_{#1}}

\def\til#1{\tilde{#1}}
\def\wtil#1{\widetilde{#1}}

\def\what#1{\widehat{#1}}

\def\sm{\setminus}

\def\E {{\mathcal E}}
\def\J {{\mathcal J}}
\def\U {{\mathcal U}}
\def\V {{\mathcal V}}
\def\z {\zeta }
\def\Harm {\text {\rm Harm}\, }
\def\grad {\nabla }
\def\dexh {\{ M_k \} _{k=0}^{\infty } }
\def\sing#1{#1_{\text {\rm sing}}}
\def\reg#1{#1_{\text {\rm reg}}}

\def\setof#1#2{\left\{\,#1\mid #2\,\right\} }
\def\set#1{{\left\{#1\right\}}}

\def\unif{_{\text{\rm unif}}}

\def\embdim{{\text{\rm embdim}}\,}

\def\dist{{\text{\rm dist}}\,}
\def\distg#1{{\text{\rm dist}_{#1}}}

\def\holecl {M\setminus \overline M_0}
\def\hole {M\setminus M_0}

\def\nd{\frac {\partial }{\partial\nu } }
\def\ndof#1{\frac {\partial#1}{\partial\nu } }

\def\pdof#1#2{\frac {\partial#1}{\partial#2}}

\def\cinf{\ensuremath{C^{\infty}} }
\def\cinfns{\ensuremath{C^{\infty}}}

\def\diam{\text {\rm diam} \, }

\def\real{\text {\rm Re}\, }

\def\imag{\text {\rm Im}\, }

\def\supp{\text {\rm supp}\, }

\def\Vol{\text {\rm vol} \, }

\def\restrict#1{\restriction_{#1}}

\def\ztanq#1{{\mathcal T}^{(#1)}}
\def\ztan{{\mathcal T}}
\def\ztanprojq#1#2{{\Pi _{{\mathcal T}^{(#1)}#2}}}
\def\ztanproj#1{{\Pi _{{\mathcal T}#1}}}

\def\ztansmoothq#1{{\widetilde{\mathcal T}}^{(#1)}}
\def\ztansmooth{{\widetilde{\mathcal T}}}
\def\ztansmoothprojq#1#2{{\Pi _{{\widetilde{\mathcal T}}^{(#1)}#2}}}
\def\ztansmoothproj#1{{\Pi _{{\widetilde{\mathcal T}}#1}}}

\def\trlevq#1#2{{\mathfrak S} ^{({#1})}(#2)}
\def\trlev#1{{\mathfrak S} (#1)}
\def\trlevdistrq#1#2#3{[{\mathfrak S}_{#1}(#2)(#3)]_{\text{\rm distr.}}}
\def\trlevdistr#1#2{[{\mathfrak S}(#1)(#2)]_{\text{\rm distr.}}}
\def\trlevadj#1{[{\mathfrak S}(\cdot )(#1)]^*}
\def\trlevprelimq#1#2#3{{\mathfrak s} ^{({#1})}(#2,#3)}
\def\trlevprelim#1#2{{\mathfrak s} (#1,#2)}

\def\trlevcpxspaceq#1#2{\widetilde{\mathfrak S} ^{({#1})}(#2)}
\def\trlevcpxspace#1{\widetilde{\mathfrak S} (#1)}

\def\zeroq#1#2{0_{{\mathcal T}^{(#1)}#2}}
\def\zeroqat#1#2#3{0_{{\mathcal T}_{#3}^{(#1)}#2}}
\def\zero#1{0_{{\mathcal T}#1}}
\def\zeroat#1#2{0_{{\mathcal T}_{#2}#1}}

\def\zeroqsome#1#2{\hat 0_{{\mathcal T}^{(#1)}#2}}
\def\zeroqatsome#1#2#3{\hat 0_{{\mathcal T}_{#3}^{(#1)}#2}}
\def\zerosome#1{\hat 0_{{\mathcal T}#1}}
\def\zeroatsome#1#2{\hat 0_{{\mathcal T}_{#2}#1}}

\def\lenv#1{\underset{#1}{\text{\rm LEnv}}\,}
\def\uenv#1{\underset{#1}{\text{\rm UEnv}}\,}

\def\analset{{\mathfrak A}}
\def\analsetsmooth{{\widetilde{\mathfrak A}}}

\def\reffamily#1#2#3#4{{(\{ {#1}_{#2}\}_{#2\in #3},\sim_{#4})}}
\def\refmetricfamily#1#2#3#4#5#6{{(\{ (#1_{#4},#2_{#4},#3_{#4})\}_{#4\in #5},\sim_{#6})}}
\def\refmetriconfamily#1#2#3#4#5#6{(\{{\mathfrak O}(#1_{#4},#2_{#4},#3_{#4})\}_{#4\in #5},\sim_{#6})}

\def\nonneg#1{{#1}_{\geq 0}}
\def\pos#1{{#1}_{>0}}

\def\plshclass{{\mathcal {W}}}
\def\strplshclass{{\mathcal S\mathcal P}}
\def\ultraplshclass{{\mathcal U\mathcal P}}
\def\plshclassquasi{{\mathcal Q\mathcal P}}
\def\strplshclassquasi{{\mathcal Q\mathcal S\mathcal P}}
\def\ultraplshclassquasi{{\mathcal Q\mathcal U\mathcal P}}

\def\plshclassgenA{{\mathcal P}_0}
\def\strplshclassgenA{{\mathcal S\mathcal P}_0}
\def\ultraplshclassgenA{{\mathcal U\mathcal P}_0}
\def\plshclassquasigenA{{\mathcal Q\mathcal P}_0}
\def\strplshclassquasigenA{{\mathcal Q\mathcal S\mathcal P}_0}
\def\ultraplshclassquasigenA{{\mathcal Q\mathcal U\mathcal P}_0}

\def\plshclassgenB{{\mathcal P}_1}
\def\strplshclassgenB{{\mathcal S\mathcal P}_1}
\def\ultraplshclassgenB{{\mathcal U\mathcal P}_1}
\def\plshclassquasigenB{{\mathcal Q\mathcal P}_1}
\def\strplshclassquasigenB{{\mathcal Q\mathcal S\mathcal P}_1}
\def\ultraplshclassquasigenB{{\mathcal Q\mathcal U\mathcal P}_1}

\def\plshclassgenC{{\mathcal P}_2}
\def\strplshclassgenC{{\mathcal S\mathcal P}_2}
\def\ultraplshclassgenC{{\mathcal U\mathcal P}_2}
\def\plshclassquasigenC{{\mathcal Q\mathcal P}_2}
\def\strplshclassquasigenC{{\mathcal Q\mathcal S\mathcal P}_2}
\def\ultraplshclassquasigenC{{\mathcal Q\mathcal U\mathcal P}_2}

\def\on#1#2#3{{\mathfrak O}(#1,#2,#3)}
\def\onat#1#2#3#4{{\mathfrak O}_{#4}(#1,#2,#3)}
\def\onof#1{{\mathfrak O}(#1)}

\def\corner#1{{\widehat{#1}}}

\def\deltanbd#1#2{{N(#1;#2)}}
\def\deltanbdg#1#2#3{{N_{#3}(#1;#2)}}



\def\anal{analytic }
\def\analns{analytic}

\def\bdd{bounded }
\def\bddns{bounded}

\def\cpt{compact }
\def\cptns{compact}

\def\cpx{complex }
\def\cpxns{complex}

\def\cont{continuous }
\def\contns{continuous}

\def\dime{dimension }
\def\dimens{dimension }

\def\exh{exhaustion }
\def\exhns{exhaustion}

\def\fn{function }
\def\fnns{function}

\def\fns{functions }
\def\fnsns{functions}

\def\holo{holomorphic }
\def\holons{holomorphic}

\def\mero{meromorphic }
\def\merons{meromorphic}

\def\holoconvex{holomorphically convex }
\def\holoconvexns{holomorphically convex}

\def\ircomp{irreducible component }
\def\concomp{connected component }
\def\ircompns{irreducible component}
\def\concompns{connected component}
\def\ircomps{irreducible components }
\def\concomps{connected components }
\def\ircompsns{irreducible components}
\def\concompsns{connected components}

\def\irred{irreducible }
\def\irredns{irreducible}

\def\con{connected }
\def\conns{connected}

\def\comp{component }
\def\compns{component}
\def\comps{components }
\def\compsns{components}

\def\mfld{manifold }
\def\mfldns{manifold}
\def\mflds{manifolds }
\def\mfldsns{manifolds}

\def\nbd{neighborhood }
\def\nbds{neighborhoods }
\def\nbdns{neighborhood}
\def\nbdsns{neighborhoods}

\def\harm{harmonic }
\def\harmns{harmonic}
\def\plh{pluriharmonic }
\def\plhns{pluriharmonic}
\def\plsh{plurisubharmonic }
\def\plshns{plurisubharmonic}

\def\qplsh#1{$#1$-plurisubharmonic}
\def\hplsh{$(n-1)$-plurisubharmonic }
\def\hplshns{$(n-1)$-plurisubharmonic}

\def\para{parabolic }
\def\parans{parabolic}

\def\rel{relatively }
\def\relns{relatively}

\def\str{strictly }
\def\strns{strictly}

\def\strg{strongly }
\def\strgns{strongly}

\def\cvx{convex }
\def\cvxns{convex}

\def\wrt{with respect to }
\def\wrtns{with respect to}

\def\st {such that }
\def\stns {such that}

\def\hm {harmonic measure }
\def\hmns {harmonic measure}

\def\hmib {harmonic measure of the ideal boundary of }
\def\hmibns {harmonic measure of the ideal boundary of}

\def\vphi{\varphi}

\def\seq#1#2{\left\{#1_{#2}\right\} }


\def\vphi {\varphi }


\def\inv{   ^{-1}  }

\def\ssp#1{^{(#1)}}

\title[Strong $q$-convexity in uniform neighborhoods]
{Strong $q$-convexity in uniform neighborhoods of subvarieties in
coverings of complex spaces}
\author[M.~Fraboni]{Michael Fraboni}
\address{Department of Mathematics and Computer Science\\Moravian College\\Bethlehem, PA 18018\\USA}
\email{mfraboni@moravian.edu}
\author[T.~Napier]{Terrence Napier$^{*}$}
\address{Department of Mathematics\\Lehigh University\\Bethlehem, PA 18015\\USA}
\email{tjn2@lehigh.edu}
\thanks{$^{*}$Research partially
supported by NSF grant DMS0306441}



\date{April 20, 2007}

\begin{abstract}
The main result is that, for any projective \cpt \anal subset $Y$
of dimension $q>0$ in a reduced \cpx space $X$, there is a \nbd
$\Omega$ of $Y$ \stns, for any covering space
\(\Upsilon\colon\widehat X\to X\) in which \(\widehat
Y\equiv\Upsilon\inv(Y)\) has no non\cpt \con \anal subsets of pure
dimension~$q$ with only \cpt \ircompsns,  there exists a $\cinf$
\exh \fn $\vphi$ on $\widehat X$ which is strongly $q$-\cvx on
$\widehat\Omega=\Upsilon\inv(\Omega)$ outside a uniform \nbd of
the $q$-dimensional \cpt \ircomps of~$\widehat Y$.
\end{abstract}

\maketitle

\section*{Introduction}\label{introduction}

According to the main result of \cite{Fraboni-Covering cpx mfld},
for any projective \cpt \anal subset $Y$ of dimension $q>0$ in a
\cpx manifold $X$, there exists a \nbd $\Omega$ of $Y$ in \(X\)
\stns, for any covering space \(\Upsilon\colon\widehat X\to X\) in
which \(\widehat Y\equiv\Upsilon\inv(Y)\) has no \cpt \ircompsns,
there exists a $\cinf$ \exh \fn $\vphi$ on $\widehat X$ which is
strongly $q$-\cvx on $\widehat\Omega=\Upsilon\inv(\Omega)$. The
case $n=1$ was obtained in \cite{Napier-coverings of proj var} and
\cite{Coltoiu-Covering of nbd}, and a similar result was first
obtained in \cite{Coltoiu-Vajaitu n-completeness} for $Y$ a fiber
of a suitable proper \holo mapping (see also~\cite{Miyazawa}). The
main goal of this paper is a more general version in which $X$ is
only a reduced \cpx space and $\widehat Y$ may have some \cpt
\ircomps (but no infinite chains of $q$-dimensional \cpt
\ircompsns).

Let $\Omega$ be an open subset of a reduced \cpx space~$X$ and let
$q$ be a positive integer.  A \fn $\vphi$ on~$\Omega$ is
\textit{\cinf strongly $q$-\cvxns} if, for each point
$p\in\Omega$, there is a proper \holo embedding~$\Phi$ of a \nbd
$U$ of $p$ in $\Omega$ into an open set $U'\subset\C^N$ and a \fn
$\vphi'\in\cinf(U')$ \st $\vphi'\circ\Phi=\vphi$ on~$U$ and \st
the Levi-form $\lev{\vphi'}$ has at most $q-1$ nonpositive
eigenvalues at each point. We also say that $\vphi$ is of class
$\strplshclass^\infty(q)$ and we write
$\vphi\in\strplshclass^\infty(q)(\Omega)$. For $q=1$, we also say
that the \fn is \textit{\cinf \str \plshns}.

There does not seem to be a single best analogous notion of weak
$q$-convexity. One natural analogue is the following. We will say
that a \fn $\vphi$ on $\Omega$ is of class $\plshclass^\infty(q)$
if, for each point $p\in\Omega$, there is a proper \holo
embedding~$\Phi$ of a \nbd $U$ of $p$ in $\Omega$ into an open set
$U'\subset\C^N$ and a \fn $\vphi'\in\cinf(U')$ \st
$\vphi'\circ\Phi=\vphi$ on~$U$ and \st the Levi-form
$\lev{\vphi'}$ has at most $q-1$ negative eigenvalues at each
point. For $q=1$, we also say that the \fn is \textit{\cinf
\plshns}.

Since every local embedding factors through the Zariski tangent
space, it follows that, if $\vphi\in\plshclass^\infty(q)(\Omega)$
($\strplshclass^\infty(q)(\Omega)$), then, for every point
$p\in\Omega$ and every local \holo model $(U,\Phi,U')$ with $p\in
U$, there is a \nbd $V'$ of $\Phi(p)$ in $U'$ and a \fn
$\vphi'\in\plshclass^\infty(q)(V')$ (respectively,
$\strplshclass^\infty(q)(V')$) with $\vphi=\vphi'\circ\Phi$ on
$\Phi\inv(V')$.

The main result of this paper is the following (see also the
stronger version Theorem~\ref{Main thm generl version}).

\begin{thm}\label{main theorem from intro}
Let $X$ be a connected reduced complex space and let $Y$ be a
\cpt \anal subset of dimension~$q>0$ in $X$. Assume that $Y$
admits a projective embedding. Then there exists a neighborhood
$\Omega$ of $Y$ in \(X\) and a discrete subset $F$ of $\R$ \stns,
for any \con covering space \(\Upsilon\colon\widehat X\to X\) in
which \(\widehat Y\equiv\Upsilon\inv(Y)\) has no non\cpt \con
\anal subsets of pure dimension~$q$ with only \cpt \ircompsns,
there exists a $\cinf$ \exh \fn $\vphi$ on $\widehat X$ which is
of class $\plshclass^\infty(q)$ on
$\what\Omega=\Upsilon\inv(\Omega)$ and of class
$\strplshclass^\infty(q)$ on $\what\Omega\sm\vphi\inv(F)$.
\end{thm}

\begin{rmks}
1. In place of $Y$ being projective, we need only assume that $Y$
admits a nowhere dense \anal subset $S$ containing $\sing Y$ \st
$Y\setminus S$ is Stein.

2. In the proof for $X$ a $2$-dimensional \cpx manifold (and $\dim
Y=1$) in \cite{Napier-coverings of proj var}, the second author
mistakenly only wrote the statement and proof for Galois
coverings. In the $2$-dimensional case, the statement and proof
are easily modified to give the version for general coverings. In
higher dimensions, more care must be taken in dealing with the
singular set of~$Y$ (which need not be discrete).
\end{rmks}

\begin{ack}
The authors would like to thank Cezar Joita and Mohan Ramachandran
for very helpful conversations.
\end{ack}

\section{The Zariski tangent space and Hermitian metrics}\label{zariski tangent space
section}

In this section we recall notions of tangent vectors and Hermitian
metrics on \cpx spaces.  Throughout this section $X$ (or $(X,\ol
_X)$) will denote a reduced \cpx space.

\subsection*{Analytic subsets}\label{analytic subsets subsection}
Unless otherwise indicated, by an {\it \anal subset} of~$X$ we
will mean a properly embedded reduced \anal subspace~$A$; that is,
a closed subset which is locally the zero set of a collection of
\holo \fns together with structure sheaf given by the restrictions
of local \holo \fns in~$X$. We will denote the structure sheaf of
$A$ by $\ol _A$ and the ideal sheaf by $\ideal A$.

By a {\it local \holo model} (or simply a {\it local model}) in
$X$ we will mean a triple $(U,\Phi ,U')$, where $U$ is an open
subset of $X$, $U'$ is an open subset of $\C ^N$ for some $N$, and
$\Phi\colon U\to U'$ is a \holo map which maps $U$ isomorphically
onto an \anal subset $\Phi (U)$ of some open subset of~$U'$. Note
that we do {\it not} require $\Phi$ to be proper, but it will be
convenient for our purposes to also have an open set $U'$
containing $\Phi (U)$ to which to refer. If $\Phi$ is a proper
embedding, then we will call $(U,\Phi ,U')$ a {\it proper local
\holo model} (or simply a {\it proper local model}). Observe that,
for $(U,\Phi,U')$ to be a local model, we do require in general
that $\Phi$ properly embed $U$ into {\it some} open subset of $U'$
(i.e. the topology on $U$ induced by $U'$ agrees with the original
topology on $U$).

\subsection*{The Zariski tangent linear space}\label{zariski tangent
linear space subsection} For each point $p\in X$, the vector space
\[
\ztanq 1_pX\equiv (\mathfrak m_p/\mathfrak m^2_p)^*,
\]
where $\mathfrak m_p$ is the maximal ideal in $\ol_X$ at $p$, is
called the {\it Zariski tangent space at} $p$. The {\it Zariski
tangent linear space}
\[
\ztanprojq 1X\colon \ztanq 1X\equiv\bigcup_{p\in X}\ztanq 1_pX\to
X
\]
has a natural reduced \cpx \anal linear space structure with the
following local models. Given a proper local model $(U,\Phi ,U')$
in $X$ with $U'\subset\C ^N$, setting $A=\Phi (U)$ we get the
\anal subset
\[
B\equiv\setof{w\in \ztanq 1U'}{z=\ztanprojq 1{U'}(w)\in A\text{
and }df(w)=0\text{ for every }f\in(\ideal A)_z} \] of $\ztanq
1U'=T^{1,0}U'=U'\times\C^N$ and a bijection $\Phi_*\colon \ztanq
1U=\ztanq 1X\restrict U=\ztanprojq 1X\inv (U)\to B$ determined by
\[
dh(\Phi_*v)=v([h\circ\Phi-h(\Phi (p))]_p)\qquad\forall\, h\in
(\ol_{\C ^N})_{\Phi (p)}, \quad p\in U, \quad v\in\ztanq 1_pU.
\]
For each point $p\in U$, the map
$(\Phi_*)_p=\Phi_*\restrict{\ztanq 1_pX}\colon\ztanq 1_pX\to
\set{\Phi (p)}\times\C^N$ is an injective \cpx linear map. The
triple $([\ztanprojq 1X]\inv (U),\Phi_*,U'\times\C^N)$ is the
proper local \holo model in $\ztanq 1X$ corresponding to the
proper local \holo model $(U,\Phi ,U')$ in $X$. Observe that we
may identify $B$ with $\ztanq 1A$ and we get a commutative diagram
of \holo maps
\begin{center}\begin{picture}(80,50)
\put(0,40){$\ztanq 1U$} \put(47,46){$\Phi_*$}
\put(32,44){\vector(1,0){40}} \put(75,40){$\ztanq 1A$}
\put(19,0){$U$}\put(47,6){$\Phi$}\put(32,4){\vector(1,0){40}}
\put(76,0){$A$}
\put(22,37){\vector(0,-1){26}}\put(82,37){\vector(0,-1){26}}
\end{picture}
\end{center}
in which $\Phi$ and $\Phi_*$ are isomorphisms. Note that
$\ztanprojq 1X:\ztanq 1X\to X$ need {\it not} be an open mapping
if $X$ is singular.

We will denote by $TX=\bigcup_{p\in X}T_pX$ the real \anal linear
space associated to the \cpx \anal linear space $\ztanq 1X$ and by
$J\colon TX\to TX$ the corresponding \cpx structure. We may
identify $\ztanq 1X=(TX,J)$ with the $(1,0)$ part $T^{1,0}X$ of
the complexification of $TX$ under the isomorphism $v\mapsto\frac
12(v-\sqrt{-1}Jv)$.

\subsection*{Tangent mappings}\label{tangent mappings subsection}
A \holo mapping $\Psi\colon X\to Y$ of reduced \cpx spaces $X$ and
$Y$ induces a \holo mapping $\Psi_*\colon \ztanq 1X\to\ztanq 1Y$
\stns, for each $p\in X$, the restriction
$(\Psi_*)_p=\Psi_*\restrict{\ztanq 1_pX}\colon \ztanq
1_pX\to\ztanq 1_{\Psi (p)}Y$ is the linear map given by
\[
(\Psi_*)_p(v)([f]_{\Psi (p)})=v([f\circ\Psi ]_p) \qquad\forall\,
v\in\ztanq 1_pX,\, f\in\mathfrak m_{Y,\Psi (p)}.
\]
If $R$ is a \fn on $\ztanq 1Y$, then we denote the \fn
$R\circ\Psi_*$ on $\ztanq 1X$ by $\Psi^*R$ (instead of
$(\Psi_*)^*R$). If $(\Psi_*)_p$ is injective for some point $p\in
X$, then the restriction of $\Psi$ to a small \nbd $U$ of $p$ in
$X$ embeds $U$ properly into some \nbd of $\Psi (p)$ in $Y$.
Moreover, $\Psi_*$ is a (proper) embedding if and only if $\Psi$
is a (proper) embedding. Finally, given another \holo mapping
$\Phi\colon Y\to Z$ to a reduced \cpx space $Z$, we have
$(\Phi\circ\Psi )_*=\Phi_*\circ\Psi_*$.

Similarly, a \cinf mapping $\Psi\colon M\to X$ of a real \cinf
manifold $M$ into $X$ induces a \cinf tangent mapping
$\Psi_*\colon TM\to TX$. For $\Psi$ a \holo mapping, the above
mappings correspond under the isomorphism $(TX,J)\cong\ztanq 1X$.

\subsection*{Embedding dimension}\label{embedding dimension
subsection} For each point $p\in X$, there is a proper local \holo
model $\Phi\colon U\to U'\subset\ztanq 1_pX\cong\C^d$ on a \nbd
$U$ of $p$ with $\Phi (p)=0$ and $\Phi_*\left(\ztanq
1_pX\right)=\ztanq 1_0\Phi (U)=\set 0\times\C^d$. Furthermore,
every \holo embedding of a \nbd of $p$ into a \cpx manifold
factors into the composition with $\Phi$ of an embedding of a \nbd
of $0$. In other words, if $\Psi\colon V\to Y$ is a \holo
embedding of a \nbd $V\subset U$ of $p$ into a \cpx manifold $Y$
and $\Psi'$ is any \holo lifting of $\Psi$ to a \nbd $V'$ of $0$
in $\C^d$ (i.e. $\Psi'\colon V'\to Y$ is a \holo map with
$\Psi'\circ\Phi=\Psi$ on a \nbd of $p$), then $\Psi'$ embeds a
\nbd of $0\in\C^d$ into a \nbd of $\Psi (p)$ in $Y$. In
particular, $d=\dim \ztanq 1_pX=\embdim_pX$ is the (minimal)
embedding dimension of the \cpx \anal space germ at $p$ determined
by $X$.

\subsection*{The associated direct sum spaces}\label{associated direct sum
spaces subsection} For a positive integer $q$, we will denote by
$\ztanq qX$ the $q$-fold reduced fiber product space
\[
\ztanprojq qX\colon\ztanq q X =\ztanq 1 X\times_{\ztanprojq
1X}\cdots\times_{\ztanprojq 1X}\ztanq 1 X\to X.
\]
Thus, for each point $p\in X$,
\[
(\ztanprojq qX)\inv (p)=\ztanq q _pX=\ztanq 1_pX\times\cdots
\times \ztanq 1_pX =\ztanq 1_pX\oplus\cdots \oplus \ztanq 1_pX
\quad (q \text{ summands}).
\]

If $\Psi\colon X\to Y$ is a \holo mapping of reduced \cpx spaces,
then, for every positive integer $q$, we define the \holo mapping
(linear on fibers) $\Psi_*^{(q)}\colon \ztanq qX\to \ztanq qY$ by
$(v_1,\dots ,v_q)\mapsto (\Psi_*v_1,\dots ,\Psi_*v_q)$. If $R$ is
a \fn on $\ztanq qY$, then we define the \fn $\Psi^*R$ on $\ztanq
qX$ by $\Psi^*R=R\circ\Psi_*\ssp q$.

\subsection*{Hermitian metrics}\label{hermitian tensors
section} A {\it Hermitian metric} $g$ in $X$ is a \cinf \fn
$g\colon\ztanq 2X\to\C$ \st $g\restrict{\ztanq 2_pX}$ is a
Hermitian inner product for each point $p\in X$. Equivalently
(see, for example, \cite{Joita-Nap-Ram -Generic q-cvx}), for each
point $p\in X$, there is a local model $(U,\Phi,U')$ in $X$ with
$p\in U$ and a (\cinfns) Hermitian metric $g'$ in $U'$ \st
\(g\restrict{\ztanq 2U}=\Phi^*g'\). We will call $g'$ a
\textit{local representation for}~$g$.  If $\Psi\colon Y\to X$ is
a local \holo embedding, then we get a pullback Hermitian metric
$\Psi^*g$ in $Y$. We will also write $g\restrict
Y=\iota^*g=g\restrict{\ztanq 2Y}$ for an inclusion $\iota\colon
Y\subset X$.

\begin{rmks}
\noindent 1. Since every local model in a \nbd of a point $p\in X$
factors through $\ztanq 1_pX$ in a \nbd of $p$, with both factors
embeddings, it follows that a local representation $g'$ as above
exists in a \nbd of $p$ in \textit{any} local model.

\noindent 2. The Hermitian metric $g$ also determines a Hermitian
metric and an associated Riemannian metric in the isomorphic \cpx
linear space $(TX,J)$.
\end{rmks}

\subsection*{Distance in a complex space}\label{Distance subsection}

Assume that $X$ is \con and let $g$ be a Hermitian metric in $X$.
Given a piecewise \cinf curve $\gamma\colon [a,b]\to X$, we define
the \textit{length of $\gamma$ with respect to $g$} by
\[
\ell_g(\gamma)=\int_a^b|\Dot\gamma(t)|_g\, dt.
\]
Given two points, $p,q\in X$, we define the \textit{distance
between $p$ and $q$ with respect to} $g$ by
\[
\distg g(p,q)=\inf\setof{\ell_g(\gamma)}{\gamma\text{ is a
piecewise \cinf curve from }p\text{ to }q\text{ in }X}.
\]
This distance is finite since such a piecewise \cinf path
connecting two given points always exists (for example, by the
existence of a resolution of singularities). We have the following
standard fact:
\begin{prop}\label{Metric modify to complete prop}
Let $X$ be \con and let $g$ be a Hermitian metric in $X$. Then
\(\distg g(\cdot,\cdot)\) is a metric (in the sense of a distance
\fnns) in $X$ which induces the given \cpx space topology in $X$.
Furthermore, there exists a positive \cont \fn $\alpha$ on $X$
\stns, if $\beta$ is any positive \cinf \fn on $X$ with
$\beta\geq\alpha$ on the complement of some \cpt subset of $X$ and
$h=\beta\cdot g$, then, for each point $p\in X$, the \fn
$x\mapsto\distg h(x,p)$ is an \exh \fnns. In particular, $\distg
h$ is a complete metric.
\end{prop}
\begin{pf}
It is easy to see that $\distg g$ is symmetric and nonnegative on
$X\times X$ and zero on the diagonal. The triangle inequality is
also easily verified. Given a point $a\in X$, we may choose a
proper local \holo model $(U,\Phi,U')$ with $a\in U$ and a
Hermitian metric $g'$ on $U'$ with $\Phi^*g'=g$ on $U$. Given a
\con \rel \cpt \nbd $V$ of $a$ in $U$, we may choose a \con \rel
\cpt \nbd $V'$ of $\Phi(a)$ in $U'$ with $\Phi\inv(V')=V$. Let
$r=\distg{g'}(\Phi(a),U'\setminus V')>0$. Given a point $x\in X\sm
V$ and a piecewise \cinf path $\gamma$ in $X$ from $a$ to $x$,
there is some $t\in (0,1]$ with $\gamma ([0,t))\subset V$ and
$\gamma(t)\in\partial V$. Thus
\(\ell_g(\gamma)\geq\ell_{g'}\left(\Phi\left(\gamma\restrict{[0,t]}\right)\right)\geq
r\). It follows that $\distg g(a,\cdot)>0$ on $X\sm\set{a}$
(since, for any given $x\in X\sm\set{a}$, we may choose such a
\nbd $V$ not containing $x$) and $B_g(a;r)\subset V$. Thus $\distg
g$ is a metric inducing a topology which is finer than the given
topology.

For the reverse containment, let $\Psi\colon\check X\to X$ be a
resolution of singularities. Thus $\check X$ is a smooth \cpx
space with \concomps $\seq {\check X}i_{i\in I}$, $\Psi$ is a
surjective proper \holo map, the \anal set $E=\Psi\inv(\sing X)$
is nowhere dense in $\check X$, $\Psi$ maps $\check X\sm E$
isomorphically onto $\reg X$, and the distinct irreducible \comps
of $X$ are given by $X_i=\Psi(\check X_i)$ for $i\in I$. We may
also choose a Hermitian metric $h$ on $\check X$, and we may let
$\check g$ be the Hermitian metric given by $\check g=\Psi^*g+h$.
Given a point $a\in X$ and a constant $r>0$, we may choose a \nbd
$U$ of the \cpt \anal set $A=\Psi\inv(a)$ in $\check X$ \stns, for
each $i\in I$, $U\cap\check X_i$ is contained in the $r$-\nbd of
$A\cap\check X_i$ with respect to the metric $\distg{\check
g\restrict{\check X_i}}$ (in particular, $U\cap\check
X_i=\emptyset$ if $A\cap\check X_i=\emptyset$). The set $\Psi(U)$
contains a \nbd $V$ of $a$ in $X$ and, given a point $\check x\in
U$, we may choose a \cinf path $\check\gamma$ in $\check X_i$ for
some $i\in I$ with $\check\gamma(0)\in A$, $\check\gamma(1)=\check
x$, and $\ell_{\check g}(\check\gamma)<r$. The path
$\gamma=\Psi(\check\gamma)$ from $a$ to $x=\Psi(\check x)$ then
satisfies
\[
\ell_g(\gamma)=\int_0^1|\dot\gamma(t)|_g\,dt\leq\int_0^1|\dot{\check\gamma}(t)|_{\check
g}\,dt<r.
\]
Thus $a\in V\subset\Psi(U)\subset B_g(a;r)$ and hence the metric
topology and the given topology are equal.

Finally, for the construction of the \fn $\alpha$, we may assume
without loss of generality that $X$ is non\cptns. We may choose a
sequence of domains $\seq\Omega\nu_{\nu=0}^\infty$ in $X$ \st
$X=\bigcup_{\nu=0}^\infty\Omega_\nu$ and
$\Omega_{\nu-1}\Subset\Omega_\nu$ for each $\nu=1,2,3,\dots$.
There then exists a \cont \fn $\alpha\colon X\to (1,\infty)$ \st
$\alpha>\left[\distg g(\Omega_{\nu-1},X\sm\Omega_\nu)\right]^{-2}$
on $X\sm\Omega_{\nu-1}$ for each $\nu=1,2,3,\dots$. Suppose $p\in
X$, $\beta$ is a positive \cinf \fn on $X$ with $\beta\geq\alpha$
on the complement of some \cpt set $K\subset X$, and $h=\beta\cdot
g$. We may fix $\mu\in\N$ with $\set{p}\cup
K\subset\Omega_{\mu-1}$. Suppose $R>0$ and $\nu\geq R+\mu$. If
$\gamma$ is a piecewise \cinf path from $p$ to a point $x\in
X\sm\Omega_\nu$, then we have numbers
\[
0<s_\mu<t_\mu\leq s_{\mu+1}<t_{\mu+1}\leq
s_{\mu+2}<t_{\mu+2}\leq\cdots\leq s_\nu<t_\nu\leq 1
\]
\stns, for $j=\mu,\dots,\nu$, we have
$\gamma((s_j,t_j))\subset\Omega_j\sm\Omega_{j-1}$,
$\gamma(s_j)\in\partial\Omega_{j-1}$, and
$\gamma(t_j)\in\partial\Omega_j$. Hence
\[
\ell_h(\gamma)\geq\sum_{j=\mu}^\nu\int_{s_j}^{t_j}|\dot\gamma(t)|_h\,dt\geq
\sum_{j=\mu}^\nu\left[\distg
g(\Omega_{j-1},X\sm\Omega_j)\right]^{-1}\int_{s_j}^{t_j}|\dot\gamma(t)|_g\,dt\geq\nu-\mu+1>R.
\]
Thus $B_h(p;R)\subset\Omega_\nu\Subset X$ and hence the \fn
$x\mapsto \distg h(x,p)$ exhausts $X$.
\end{pf}

To close this section, we record two facts for later use.
\begin{lem}\label{complete on covering lem}
Suppose $X$ is \conns, $g$ is a Hermitian metric on $X$, and
$x\mapsto\dist_g(x,p)$ is an \exh \fn on $X$ for some (hence, for
each) point $p\in X$. Then we have the following:
\begin{enumerate}
\item[(a)] For each point $p\in X$ and each constant $R>0$, the
set
$$
K(p,R)\equiv \setof{[\alpha ]\in \pi _1(X,p)}{\alpha \text{ is a
piecewise }\cinf \text{ loop in }X \text{ of length }<R}
$$
is finite.

\item[(b)] For every \con covering space $\Upsilon\colon\what X\to
X$, $x\mapsto\dist_{\hat g}(x,p)$ is an \exh \fn for each point
$p\in\what X$, where $\hat g=\Upsilon^*g$. In particular,
$\dist_{\hat g}(\cdot,\cdot)$ is a complete metric on~$\what X$.
\end{enumerate}
\end{lem}
\begin{pf}
For the proof of (a), we fix a point $p\in X$ and number $r>0$ so
small that, for each point $a$ in the \cpt set
$D\equiv\overline{B_g(p;R)}$, the ball $B_g(a;3r)$ is contained in
some contractible open set. We may also choose points
$p=p_1,p_2,\dots ,p_k\in D$ \st the balls $B_1=B_g(p_1;r),\dots,
B_k=B_g(p_k;r)$ form a covering for $D$, a Lebesgue number
$\delta>0$ for this covering with $\delta<r$, and a positive
integer $m$ \st $R/m<\delta$. For each pair of indices $i,j$, we
may choose a piecewise \cinf path $\lambda_{ij}=\lambda_{ji}\inv$
from $p_i$ to $p_j$ \st $\ell_g(\lambda_{ij})<\dist_g(p_i,p_j)+r$.
Now any piecewise $\cinf$ loop~$\alpha$ of length~$<R$ based at
$p$ is homotopic to a loop $\alpha _1 *\alpha _2*\cdots *\alpha
_m$ in $D$; where, for each $\nu=1,\dots,m$, $\alpha _\nu$ is a
piecewise \cinf path of length $<\delta$ and is, therefore,
contained in $B_{i_\nu}$ for some index $i_\nu$. We may choose
$i_1=i_m=1$. For each $\nu=2,\dots,m$, we may choose a piecewise
\cinf path $\rho_\nu$ of length~$<r$ from
$\alpha_{\nu-1}(1)=\alpha_\nu(0)$ to~$p_{i_{\nu-1}}$. Thus \(
\alpha\sim\gamma_1*\cdots*\gamma_m\), where
$\gamma_1=\alpha_1*\rho_2$, $\gamma_m=\rho_m\inv*\alpha_m$, and,
for $\nu=2,\dots,m-1$,
$\gamma_\nu=\rho_\nu\inv*\alpha_\nu*\rho_{\nu+1}$. For each
$\nu=2,\dots,m$, we have $\ell_g(\gamma_\nu)<3r$. We also have
$\alpha_\nu(0)\in B_{i_{\nu-1}}\cap B_{i_\nu}$, so
$\ell_g\left(\lambda_{i_{\nu-1}i_\nu}\right)<\dist_g(p_{i_{\nu-1}},p_{i_\nu})+r<3r$.
Therefore, $\gamma_\nu$ is homotopic to
$\lambda_{i_{\nu-1}i_\nu}$. Moreover, $\gamma_1$ is homotopic to
the trivial loop at~$p_1=p$. Thus $\alpha$ is homotopic to the
loop
$\lambda_{i_1i_2}*\lambda_{i_2i_3}*\cdots*\lambda_{i_{m-1}i_m}$
and the claim follows.

For the proof of~(b), suppose $\Upsilon\colon\what X\to X$ is a
\con covering space, $\hat g=\Upsilon^*g$, $p\in X$, $R>0$, and
$\seq x\nu$ is a sequence in $\what X$ \st $\dist_{\hat
g}(x_\nu,p)<R$ for each~$\nu$. We must show that $\seq x\nu$
admits a convergent subsequence. Since
$\set{\Upsilon(x_\nu)}\subset B_g(p;R)\Subset X$, we may assume
that $\set{\Upsilon(x_\nu)}$ converges to some point~$a\in X$.
Fixing a contractible \nbd $U$ of~$a$ in $X$ and a sufficiently
small constant $\epsilon>0$, we may also assume that
$\{\Upsilon(x_\nu)\}\subset B_g(a;\epsilon)\Subset U$. Therefore,
by replacing $R$~with~$R+\epsilon$ and each point $x_\nu$ with the
unique point in $\Upsilon\inv(a)$ lying in the same \concomp of
$\Upsilon\inv(U)$ as~$x_\nu$, we may assume that
$\Upsilon(x_\nu)=a$ for each~$\nu$. For each $\nu$, there exists a
piecewise \cinf path $\gamma_\nu$ of length~$<R$ from $p$ to
$x_\nu$. Each of the homotopy classes
$\left[\Upsilon(\gamma_1\inv*\gamma_\nu)\right]$ in $\pi_1(X,a)$
lies in the finite set $K(a,2R)$, so the collection of terminal
points $\setof{x_\nu}{\nu=1,2,3,\dots}$ of the liftings
$\set{\gamma_1\inv*\gamma_\nu}$ must be a finite set. The claim
now follows.
\end{pf}

\begin{lem}\label{Uniform dist between comps in cover lem}
Let $(X,g)$ be a \con reduced Hermitian \cpx space, let $Y$ be a
\cpt \anal subset of~$X$, and let $\epsilon>0$. Then there exists
a constant $\delta>0$ \stns, for every \con covering space
$\Upsilon\colon\what X\to X$ and for every pair of \ircomps
$A$~and~$B$ of $\what Y=\Upsilon\inv(Y)$, we have $\dist_{\hat
g}(A\sm \deltanbd{A\cap B}\epsilon,B)>\delta$; where $\hat
g=\Upsilon^*g$. In particular, for $A$ and $B$ disjoint, we have
$\dist_{\hat g}(A,B)>\delta$.
\end{lem}
\begin{rmk}
If $A\subset\deltanbd{A\cap B}\epsilon$, then we take $\dist_{\hat
g}(A\sm\deltanbd{A\cap B}\epsilon,B)=\dist_{\hat
g}(\emptyset,B)=\infty$.
\end{rmk}
\begin{pf}
We may choose open sets $\seq Vj_{j=1}^m$ in $X$ and contractible
open sets $\seq Uj_{j=1}^m$ in $X$ \st
$Y\subset\bigcup_{j=1}^mV_j$ and \stns, for each $j=1,\dots,m$, we
have $V_j\Subset U_j$. We may fix a Lebesgue number~$\eta>0$ for
the covering $\seq Vj_{j=1}^m$ of $Y$ with respect to the metric
$\dist_g(\cdot,\cdot)$. In other words, for each point $p\in Y$,
we have $B_g(p;\eta)\subset V_j$ for some~$j$. We may also choose
$\eta$ so that $\eta<\epsilon$. We may then choose $\delta>0$ so
that $2\delta<\eta$ and so that, for each $j=1,\dots,m$, we have
$2\delta<\dist_g(\overline V_j\cap A\sm\deltanbd{A\cap B}\eta,B)$
for every pair of \ircomps $A$ and $B$ of $Y\cap U_j$.

Now suppose that $\Upsilon\colon\what X\to X$ is a \con covering
space, $\hat g=\Upsilon^*g$, $A$ and $B$ are \ircomps of~$\what
Y=\Upsilon\inv(Y)$, and $\hat a\in A$ and $\hat b\in B$ with
$\dist_{\hat g}(\hat a,\hat b)<2\delta$. Setting $a=\Upsilon(\hat
a)$ and $b=\Upsilon(\hat b)$, we then have $b\in
B_g(a;2\delta)\subset B_g(a;\eta)\subset V_j$ for some~$j$. Hence
\[
\hat b\in B_{\hat g}(\hat a;2\delta)\subset B_{\hat g}(\hat
a;\eta)\subset V\subset U,
\]
where $U$ is the \concomp of $\Upsilon\inv(U_j)$ containing~$\hat
a$ (which $\Upsilon$ maps isomorphically onto~$U_j$) and
$V=\Upsilon\inv(V_j)\cap U$. We have $\hat a\in\what A_0$ and
$\hat b\in\what B_0$ for some \ircomps $\what A_0$ of $A\cap U$
and $\what B_0$ of $B\cap U$ (and, therefore, of $\what Y\cap U$).
Since $\Upsilon$ maps $U$ isomorphically onto $U_j$, the sets
$A_0=\Upsilon(\what A_0)$ and $B_0=\Upsilon(\what B_0)$ are
\ircomps of $Y\cap U_j$. Since $\dist_g(a,b)<2\delta$, the choice
of $\delta$ implies that $a\in \overline V_j\cap
A_0\cap\deltanbd{A_0\cap B_0}\eta$. Hence we may choose a
piecewise \cinf path $\lambda$ of length less than~$\eta$ in $X$
with $\lambda(0)=a$ and $\lambda(1)\in A_0\cap B_0$. The lifting
$\hat\lambda$ with $\hat\lambda(0)=\hat a$ lies entirely in $U$
and, therefore, $\hat\lambda(1)\in\what A_0\cap\what B_0\subset
A\cap B$. Hence $\dist_{\hat g}(\hat a,A\cap B)<\eta<\epsilon$ and
it follows that $\dist_{\hat g}(A\sm\deltanbd{A\cap
B}\epsilon,B)\geq 2\delta>\delta$.
\end{pf}

\section{Positive sums of eigenvalues of the Levi form}\label{Pos sums of ev of the Levi sect}

Throughout this section $X$ will denote a reduced \cpx space, $g$
will denote a Hermitian metric in $X$, and $q$ will denote a
positive integer.
\begin{defn}\label{plsh wrt g,q def}
Let $\vphi$ be a real-valued \fn on an open subset $\Omega$ of
$X$. We will say that $\vphi$ is \textit{of class
$\plshclass^\infty(g,q)$} (\textit{of class
$\strplshclass^\infty(g,q)$}) and write
$\vphi\in\plshclass^\infty(g,q)(\Omega)$ (respectively,
$\vphi\in\strplshclass^\infty(g,q)(\Omega)$) if, for every \anal
subset $Y$ of an open subset of $\Omega$, every point $p\in Y$,
every local \holo model $(U,\Phi,U')$ in $Y$ with $p\in U$, and
every Hermitian metric $g'$ in $U'$ with $\Phi^*g'=g\restrict Y$
on $U$, there exists a \cinf \fn $\vphi'$ on a \nbd $V'$ of
$\Phi(p)$ in $U'$ \st $\vphi'\circ\Phi=\vphi$ on
$V=\Phi\inv(V')\subset U$ and \stns, for each point $z\in V'$, the
trace of the restriction of the Levi-form \(\lev{\vphi'}\) to any
$q$-dimensional subspace of $\ztanq 1_zV'$ with respect $g'$ is
nonnegative (respectively, positive); that is, for any
$g'$-orthonormal collection of vectors $e_1,\dots,e_q\in\ztanq
1_zV'$, we have
\[
\sum_{i=1}^q\lev{\vphi'}(e_i,e_i)\geq 0 \quad\text{(respectively,
}>0\text{)}.
\]
\end{defn}

Clearly, $\plshclass^\infty(g,q)\subset\plshclass^\infty(q)$ and
$\strplshclass^\infty(g,q)\subset\strplshclass^\infty(q)$. For
manifolds, these classes of \fns were first introduced by Grauert
and Riemenschneider~\cite{Grauert-Riemenschneider} and have been
applied by others in many different contexts. The definition is
stated in the above (rather cumbersome) form in order to guarantee
the above inclusions as well as to guarantee invariance under sums
(by independence of the local representation) and restrictions to
\anal subsets. One weakness in the definition of
$\plshclass^\infty(g,q)$ is that the local extension~$\vphi'$ is
not assumed to be of class~$\plshclass^\infty(g',q)$; i.e.~it is
not assumed that $\vphi'$ admits further extensions with Levi form
satisfying the trace condition with respect to other local
representations of~$g'$. One could build the existence of such
extensions into the definition (in fact, the \fns produced in this
paper will actually have such properties), but the above notion
suffices for our purposes. For $\strplshclass^\infty(g,q)$, the
situation is much better. In fact, the following proposition
provides an equivalent notion which is more easily checked (see,
for example, \cite{Joita-Nap-Ram -Generic q-cvx} as well as
\cite{Richberg}, \cite{Demailly-Cohomology of q-convex spaces},
\cite{Coltoiu-Complete locally pluripolar sets}, \cite{Joita
Traces of cvx domains}):

\begin{prop}\label{str q-plsh equiv def prop}
A \fn $\vphi$ on an open subset $\Omega$ of $X$ is of class
$\strplshclass^\infty(g,q)$ if and only if, for every point
$p\in\Omega$, there exist a local \holo model $(U,\Phi,U')$ in $X$
with $p\in U\subset\Omega$ and a \cinf \fn $\vphi'$ on $U'$ \st
$\vphi'\circ\Phi=\vphi$ on $U$ and \stns, for each point $x\in U$,
the trace of the restriction of the pullback of the Levi-form
$\Phi^*\lev{\vphi'}$ to any $q$-dimensional subspace of $\ztanq
1_xU$ with respect to $g$ is positive.
\end{prop}

The analogous class in which the trace of the restriction of
$\Phi^*\lev{\vphi'}$ is only assumed to be nonnegative is
considered in \cite{Joita-Nap-Ram -Generic q-cvx}.  The authors do
not know whether or not this gives the same class
$\plshclass^\infty(g,q)$.

One also has the following global extension theorem of
Richberg~\cite{Richberg} (cf. Fritzsche~\cite{Fritzsche},
Demailly~\cite{Demailly-Cohomology of q-convex spaces}, Col\c
toiu~\cite{Coltoiu-Complete locally pluripolar sets},
Joita~\cite{Joita Traces of cvx domains}, and \cite{Joita-Nap-Ram
-Generic q-cvx}):
\begin{thm}\label{Global extension of str q-plsh thm}
If $Y$ is an \anal subset of $X$ and
$\vphi\in\strplshclass^\infty(g\restrict Y,q)(Y)$, then there
exists a \fn $\psi$ of class $\strplshclass^\infty(g,q)$ on a \nbd
of $Y$ in $X$ \st $\psi=\vphi$ on $Y$.
\end{thm}

Combining this with modifications of the arguments of Greene and
Wu~\cite{Greene-Wu embeddings}, Demailly~\cite{Demailly-Cohomology
of q-convex spaces}, and Ohsawa~\cite{Ohsawa Completeness 1984},
one gets the following (see \cite{Joita-Nap-Ram -Generic q-cvx}):
\begin{thm}\label{Qplsh on nbd of lower dimension or noncpt q dim
thm} Let $Y$ be a (properly embedded) \anal subset of $X$ of
dimension~$\leq q$ with no \cpt \ircomps of dimension~$q$. Then
there exists a \cinf \exh \fn $\vphi$ on $X$ which is of class
$\strplshclass^\infty(g,q)$ on a \nbd of $Y$.
\end{thm}

It will also be convenient to have the following easy lemma:
\begin{lem}\label{convex fn comp lemma}
Let $\chi$ be a \cinf \fn on an interval $(a,b)$ in $\R$ with
$\chi'\geq 0$ and $\chi''\geq 0$ and let $\vphi$ be a \fn of class
$\strplshclass^\infty(g,q)$ on an open set $\Omega\subset X$ with
values in $(a,b)$. Then the \fn $\psi\equiv\chi(\vphi)$ has the
following properties:
\begin{enumerate}
\item[(a)] We have $\psi\in\plshclass^\infty(g,q)(\Omega)$.

\item[(b)] If $\chi'>0$, then
$\psi\in\strplshclass^\infty(g,q)(\Omega)$.

\item[(c)] For any \cinf \fn $\alpha$ with \cpt support in
$\Omega$, there is a constant $C=C(\vphi,\alpha)>0$ \st
$\psi+\alpha$ will be of class $\strplshclass^\infty(g,q)$ on
$\Omega$ for any choice of $\chi$ with $\chi'>0$ (and $\chi''\geq
0$) on $(a,b)$ and $\chi'>C$ on $\vphi(\supp\alpha)$.

\end{enumerate}
\end{lem}
\begin{rmk}
Most of the $\plshclass^\infty(g,q)$ \fns to be constructed in
this paper will locally be expressible as sums of functions of the
form $\psi=\chi(\vphi)$ as in the above lemma.
\end{rmk}

\section{$Q$-convexity in a covering of a neighborhood of a Stein manifold}\label{q-cvx on nbd off sing section}

In this section, we produce, in a covering, a \fn on a uniform
\nbd of the lifting of an embedded Stein manifold which is \strg
$q$-\cvx near the lifting of a large \cpt set and equal to~$0$
elsewhere.  More precisely, we prove the following:
\begin{prop}\label{q-cvx on nbd off sing set prop}
Let $(X,g)$ be a \con reduced Hermitian \cpx space of bounded
local embedding dimension; let $Y$ be a \con Stein manifold of
dimension~$q>0$ which is properly embedded in $X$; and let $Q$,
$\Omega_1$, and $\Omega_2$ be open subsets of $X$ \st
\[
X\Supset\Omega_1\Supset\Omega_2\Supset Q;
\]
$Q\cap Y$ and $\Omega_j\cap Y$ for $j=1,2$ are nonempty, \conns,
and \rel \cpt in $Y$; $\overline {Q\cap Y}=\overline Q\cap Y$; and
$\overline {\Omega_j\cap Y}=\overline \Omega_j\cap Y$ for $j=1,2$.
Then there exists a \con \nbd $\Theta_1$ of $Y$ in $X$ and a
nonnegative \cinf \fn $\alpha$ on $\Theta_1$ satisfying the
following.
\begin{enumerate}
\item[(a)] We have
\begin{enumerate}
    \item[(i)] On $\Theta_1\sm\Omega_1$, $\alpha\equiv 0$;
    \item[(ii)] On $\Theta_1\cap\Omega_2$, $\alpha>0$;
    \item[(iii)] On some \nbd of $\Theta_1\sm Q$, $\alpha$ is of
    class $\plshclass^\infty(g,q)$; and
    \item[(iv)] For $B=\setof{x\in\Theta_1}{\alpha(x)>0}\supset\Theta_1\cap\Omega_2$,
    $\alpha$ is of class $\strplshclass^\infty(g,q)$ on a \nbd of $B\sm Q$.
\end{enumerate}
\item[(b)] Suppose $\Omega_3$ and $\Omega_4$ are open subsets of
$X$ \st
\[
\Omega_2\Supset\Omega_3\Supset\Omega_4\Supset Q,
\]
and, for $j=3,4$, $\Omega_j\cap Y$ is \con and
$\overline{\Omega_j\cap Y}=\overline\Omega_j\cap Y$. Suppose also
that $\seq\omega\nu_{\nu\in N}$ is a family of real-valued \cinf
\fns with \cpt support in $B$. Then, for some \con \nbd $\Theta_2$
of $Y$ in $\Theta_1$, for any family of sufficiently large
positive constants $\seq R\nu_{\nu\in N}$, for every \con infinite
covering space \linebreak $\Upsilon:\what X\to X$ in which $\what
Y=\Upsilon\inv(Y)$ and $\Upsilon\inv(\Omega_4\cap Y)$ are \conns,
and for every positive \cont \fn $\theta$ on $\what X$, there is a
nonnegative \cinf \fn $\alpha_0$ on
$\what\Theta_2=\Upsilon\inv(\Theta_2)$ \stns, if
$\what\Omega_j=\Upsilon\inv(\Omega_j)$ for $j=1,2,3,4$,
$\hat\alpha=\alpha\circ\Upsilon$, and $\what B=\Upsilon\inv(B)$,
then
\begin{enumerate}
    \item[(i)] On $\what\Theta_2$, $\alpha_0\geq\hat\alpha$;
    \item[(ii)] On $\what\Theta_2\sm\what\Omega_3$, $\alpha_0=\hat\alpha$;
    \item[(iii)] For each point $\nu\in N$, the \fn
    $R_\nu\cdot\alpha_0+\omega_\nu\circ\Upsilon$ will be
    \cinf \strg $q$-\cvx on $\what B\cap\what\Theta_2$; and
    \item[(iv)] On $\what\Theta_2\cap\what\Omega_4$,
    $\alpha_0>\theta$.
\end{enumerate}
\end{enumerate}
\end{prop}

\begin{lem}\label{Extension special lemma}
Let $(X,g)$ be a Hermitian \cpx manifold, let $Z$ be a (properly
embedded) \cpx submanifold, and let $\alpha$ be a real-valued
$\cinf$ \fn on $X$ \st $\alpha\restrict Z$ is of class
$\strplshclass^\infty(g\restrict Z,q)$ on $Z$. Suppose $\rho$ is a
\cinf \fn on $X$ \stns, for each point $p\in Z$, we have
$\rho(p)=0$, $(d\rho)_p=0$, and $\lev\rho(v,v)>0$ for each tangent
vector $v\in T^{1,0}_pX\sm T^{1,0}Z$ (in particular, $\lev\rho
(v,v)\geq 0$ for each $v\in T^{1,0}_pX$). Then there exists a
positive \cont \fn $\lambda_0$ on $Z$ \stns, for every \fn
$\lambda\in\cinf(X)$ with $\lambda>\lambda_0$ on $Z$, the \fn
$\beta\equiv\alpha+\lambda\cdot\rho$ is of class
$\strplshclass^\infty(g,q)$ on a \nbd of $Z$ in $X$ (depending on
the choice of $\lambda$).
\end{lem}
\begin{pf}
Suppose $\lambda$ is a positive \cinf \fn on $X$ and
$\beta\equiv\alpha+\lambda\cdot\rho$. For each point $p\in Z$ and
each tangent vector $v\in T^{1,0}_pX$, we have
\[
\lev\beta(v,v)=\lev\alpha(v,v)+\lambda(p)\cdot\lev\rho(v,v)\geq
\lev\alpha(v,v).
\]
Let $E_X$ and $E_Z$ be the set of orthonormal $q$-frames in
$\ztanq qX$ and $\ztanq qZ$, respectively. By the definition of
$\strplshclass^\infty(g,q)$ and continuity, there is a \nbd $W$ of
$E_Z$ in $\ztanq qX$ \st
\[
\sum_{j=1}^q\lev\alpha(v_j,v_j)>0 \qquad\forall\,
v=(v_1,\dots,v_q)\in W.
\]
On the other hand, for each $(e_1,\dots,e_q)\in (E_X\setminus
W)\cap[\ztanprojq qX]\inv (Z)$, we have
\[
\sum_{j=1}^q\lev\rho(e_j,e_j)>0.
\]
Therefore, if $\lambda$ grows sufficiently quickly at infinity on
$Z$, then we will have
\[
\sum_{j=1}^q\lev\beta(e_j,e_j)>0\qquad\forall\, (e_1,\dots
,e_q)\in E_X\cap[\ztanprojq qX]\inv (Z),
\]
and it follows that $\beta$ will be of class
$\strplshclass^\infty(g,q)$ on a \nbd of $Z$.
\end{pf}

\begin{lem}\label{Existence of good fn on Stein lemma}
Let $X$ be a Stein manifold and let $Z$ be a (properly embedded)
\cpx submanifold.  Then there exists a nonnegative \cinf \plsh
\fnns~$\rho$ on $X$ \st $\rho(p)=0$ and $(d\rho)_p=0$ for each
point $p\in Z$ and $\lev\rho (v,v)>0$ for each nonzero tangent
vector $v\in T^{1,0}X\setminus T^{1,0}Z$.
\end{lem}
\begin{pf}
For each point $p\in X$, Cartan's Theorem~A provides \holo \fns
\[
f_p^{(1)},\dots,f_p^{(m_p)}\in H^0(X,\cal I_{\set{p}\cup Z})
\]
generating the ideal sheaf $\cal I_{\set{p}\cup Z}$ at each point
in a \nbd $U_p$ of $p$ in $X$. Setting
$f_p=\left(f_p^{(1)},\dots,f_p^{(m_p)}\right)$, we see that, for
any nonzero vector $v\in T^{1,0}U_p$, we have $df_p(v)=0$ if and
only if $v\in T^{1,0}Z$. Forming a countable cover
$\set{U_{p_\nu}}$ of $X$ by such sets and choosing a sequence of
positive numbers $\seq\epsilon\nu$ converging to~$0$ sufficiently
fast, we get the \fn $\rho\equiv\sum\epsilon_\nu |f_{p_\nu}|^2$
with the required properties.
\end{pf}

As in \cite{Fraboni-Covering cpx mfld}, we will also apply the
following theorem of Demailly~\cite{Demailly-Cohomology of
q-convex spaces} (see Theorem~1.13 of \cite{NR-Exhausting
subsolutions}):
\begin{thm}\label{Exh subsolutions Demailly NR thm}
Let $D$ be a \con closed non\cpt subset of a Hermitian manifold
$(X,g)$, let $U$ be a \con \nbd of $D$ in $X$, and let $\theta$ be
a positive \cont \fn on $X$. Then there exists a \cinf sub\harm
\fn $\alpha$ on $X$ \st $\alpha\equiv 0$ on $X\sm U$, $\alpha>0$
and $\lap\alpha >0$ on $U$, and $\alpha>\theta$ and
$\lap\alpha>\theta$ on $D$.
\end{thm}
\begin{pf*}{Proof of Proposition~\ref{q-cvx on nbd off sing set
prop}} We first prove the proposition for $X$ a \con open subset
of $\C^n$ and $Y$ closed in $\C^n$. Let $\pi\colon\cal
N=[T^{1,0}\C^n\restrict Y]/T^{1,0}Y\to Y$ be the normal bundle
with the Hermitian metric $h$ induced by the Euclidean metric.
Then there exists a positive \cont \fn $\sigma$ on $Y$ and a
biholomorphism $\Psi\colon N\to M$ of the \nbd
\[
N=\setof{\xi\in\cal N}{|\xi|_h<\sigma(\pi(\xi))}
\]
of the $0$-section onto a \con \nbd $M$ of $Y$ in $X$ \st
$\Psi(0_y)=y$ for each point $y\in Y$ (see, for example,
\cite{Forster-Ramspott-Holo projection}). Setting
$\Lambda(t,x)=\Lambda_t(x)=\Psi(t\Psi\inv(x))$ for each $(t,x)\in
[0,1]\times M$, we get a \cinf strong deformation retraction
$\Lambda\colon [0,1]\times M\to M$ of $M$ onto $Y$ for which the
map $\Lambda_t\colon M\to M$ is \holo for each $t\in [0,1]$. In
particular, $\Lambda_0\colon M\to Y$ is a \holo submersion which
is equal to the identity on $Y$.

According to Lemma~\ref{Existence of good fn on Stein lemma},
there exists a nonnegative \cinf \plsh \fn $\rho$ on $\C^n$ \st
$\rho(z)=0$ and $(d\rho)_z=0$ for each point $z\in Y$ and
$\lev{\rho}(v,v)>0$ for each nonzero tangent vector $v\in
T^{1,0}\C^n\setminus T^{1,0}Y$.

Fixing a point $p\in Q\cap Y$, applying Theorem~\ref{Exh
subsolutions Demailly NR thm} in the Hermitian manifold
$Y\sm\set{p}$ to the \con relatively open set $\Omega_1\cap
Y\sm\set{p}$ and the non\cpt \con relatively closed set
$\overline\Omega_2\cap Y\sm\set{p}\subset\Omega_1\cap
Y\sm\set{p}$, and patching with a positive \cinf \fn near $p$ in
$Q$, we get a nonnegative \cinf \fn $\beta$ on $Y$ \st
\begin{enumerate}
    \item[(\ref{q-cvx on nbd off sing set prop}.1)] On $Y\setminus \Omega_1$, $\beta\equiv 0$;
    \item[(\ref{q-cvx on nbd off sing set prop}.2)] On $\Omega_1\cap Y$, $\beta>0$; and
    \item[(\ref{q-cvx on nbd off sing set prop}.3)] On $\Omega_1\cap Y\setminus Q$, $\lap_g\beta>0$.
\end{enumerate}
We may now choose a constant $\delta$ with
$0<3\delta<\min_{\overline \Omega_2\cap Y}\beta$ and open sets
$U_1$ and $Q_0$ in $X$ \st
\[
\Omega_1\Supset U_1\Supset\Omega_2, \quad Q\Supset Q_0,\quad
\beta<\delta\text{ on }Y\sm U_1, \text{ and } \lap_g\beta>0\text{
on }\Omega_1\cap Y\setminus Q_0.
\]
According to Lemma~\ref{Extension special lemma}, if we fix a
positive \fn $\lambda\in\cinf(\Omega_1)$ which is constant on
$U_1$ and sufficiently large on $\Omega_1\cap Y$, then the \fn
$\beta\circ\Lambda_0+\lambda\rho$ will be of class
$\strplshclass^\infty(g,q)$ on a \nbd of $U_2\setminus Q_0$ for
some \nbd $U_2$ of $\Omega_1\cap Y$ in $\Omega_1\cap M$.
Furthermore, choosing $U_2$ sufficiently small, we get
$\beta\circ\Lambda_0+\lambda\rho<\delta$ ($>3\delta$) on $U_2\sm
U_1$ (respectively, on $U_2\cap\Omega_2$). Therefore, choosing an
open set $U_3$ in $X$ with $Y\sm\Omega_1\subset
U_3\subset\overline U_3\subset M\setminus\overline U_1$; a \cinf
\fn $\chi\colon\R\to [0,\infty)$ \st $\chi'\geq 0$ and $\chi''\geq
0$ on $\R$, $\chi(t)\equiv 0$ for $t\leq\delta$, and
$\chi(t)=t-2\delta$ for $t\geq 3\delta$; and a \con \nbd
$\Theta_1$ of $Y$ in $U_3\cup U_2$, we may define a \cinf \fn
$\alpha$ on $\Theta_1$ by setting $\alpha=0$ on $\Theta_1\sm U_2$
and $\alpha=\chi(\beta\circ\Lambda_0+\lambda\rho)$ on
$\Theta_1\cap U_2$. It is now easy to verify that $\alpha$ has the
properties described in (a) (using Lemma~\ref{convex fn comp
lemma} to get~(iii)).

For the construction of the \nbd $\Theta_2$ as in (b), we choose
open sets $V_1$, $V_2$, and $V_3$ in $X$ \st \(\Omega_3\Supset
V_1\Supset V_2\Supset V_3\Supset\Omega_4\) and \stns, for
$j=1,2,3$, $V_j\cap Y$ is \con and $\overline{V_j\cap Y}=\overline
V_j\cap Y$. If $\sigma_0<\sigma$ is a sufficiently small positive
\cont \fn $Y$, then the \con open set
$\Theta_2\equiv\Psi(\setof{\xi\in\cal
N}{|\xi|_h<\sigma_0(\pi(\xi))})$, will satisfy
\[
Y\subset\Theta_2\subset\overline\Theta_2\subset\Theta_1,\quad
\Lambda_0\inv(V_2\cap Y)\cap\Theta_2\subset V_1,\text{ and }
\Lambda_0(\Omega_4\cap\Theta_2)\subset V_3\cap Y.
\]

We may also choose a nonnegative \cinf \fn $\gamma$ on $\C^n$ \st
$\supp\gamma\subset\Omega_3\cap\Theta_1$, $\gamma $ is \cinf \str
\plsh on a \nbd of $\overline V_1\cap\overline\Theta_2$, and
$\alpha+\gamma$ is of class $\strplshclass^\infty(g,q)$ on a \nbd
of $B\cap\overline\Theta_2\setminus Q$; where
$B=\setof{x\in\Theta_1}{\alpha(x)>0}\supset\Theta_1\cap\Omega_2$.

Suppose now that $\Upsilon:\what X\to X$ is a \con infinite
covering space in which $\what Y=\Upsilon\inv(Y)$ and
$\Upsilon\inv(\Omega_4\cap Y)$ are \conns.  Let $\theta$ be a
positive \cont \fn on $\what X$, let $\what M=\Upsilon\inv (M)$,
let $\what Q=\Upsilon\inv(Q)$,  let
$\what\Theta_j=\Upsilon\inv(\Theta_j)$ for $j=1,2$, let
$\what\Omega_j=\Upsilon\inv(\Omega_j)$ for $j=1,2,3,4$, let $\what
V_j=\Upsilon\inv(V_j)$ for $j=1,2,3$, let
$\hat\alpha=\alpha\circ\Upsilon$, let
$\hat\beta=\beta\circ\Upsilon\restrict{\what Y}$, let
$\hat\rho=\rho\circ\Upsilon$, let
$\hat\lambda=\lambda\circ\Upsilon$, and let $\hat g=\Upsilon^*g$.
Observe that the map $\Lambda$ lifts to a unique \cinf strong
deformation retraction $\hat\Lambda$ of $\what M$ onto $\what Y$.
The restriction of the corresponding \holo submersion
$\hat\Lambda_0\colon\what M\to\what Y$ to
$\overline{\what\Theta_2}$ is a proper map and hence we may choose
a positive \cont \fn $\theta_0$ on $\what Y$ \st
$\theta_0\circ\hat\Lambda>\theta$ on $\what\Theta_2$.

Applying Theorem~\ref{Exh subsolutions Demailly NR thm} in the
Hermitian manifold $(\what Y,\hat g\restrict{\what Y})$ to the
\con relatively open set $\what V_2\cap\what Y$ and the non\cpt
\con closed set $\overline{\what V_3}\cap\what Y$, we get a
nonnegative \cinf \fn $\beta_0$ on $\what Y$ \st
\begin{enumerate}
    \item[(\ref{q-cvx on nbd off sing set prop}.4)] On $\what Y\sm\what V_2$, $\beta_0\equiv 0$;
    \item[(\ref{q-cvx on nbd off sing set prop}.5)] On $\what V_2\cap\what Y$, $\beta_0>0$ and $\lap_{\hat g}\beta_0>0$; and
    \item[(\ref{q-cvx on nbd off sing set prop}.6)] On $\what V_3\cap\what Y$, $\beta_0>\theta_0$ and
    $\lap_{\hat g}(\hat \beta+\beta_0)>0$.
\end{enumerate}
In particular, $\lap_{\hat g}(\hat \beta+\beta_0)>0$ on
$\what\Omega_1\cap\what Y$. Setting $\what B=\Upsilon\inv(B)$,
$\hat\gamma=\gamma\circ\Upsilon$, and
$\alpha_1=(\hat\alpha+\beta_0\circ\hat\Lambda_0)\restrict{\what\Theta_2}$,
we will show that the \cinf \fn $\alpha_0=\alpha_1+\hat\gamma$ has
the required properties.

Clearly, we have $\alpha_0\geq\hat\alpha$ on $\what\Theta_2$.
Since $\alpha_1=\hat\alpha$ on
$\what\Theta_2\sm\hat\Lambda_0\inv(\what V_2\cap\what
Y)\supset\what\Theta_2\sm\what
V_1\supset\what\Theta_2\sm\what\Omega_3$ and $\supp\gamma\subset
\Omega_3$, we have $\alpha_0=\hat\alpha$ on
$\what\Theta_2\sm\what\Omega_3$. On the set
$\hat\Lambda_0\inv(\what V_3\cap\what
Y)\cap\what\Theta_2\supset\what\Omega_4\cap\what\Theta_2$, we have
\(\alpha_0\geq\beta_0\circ\hat\Lambda_0>\theta_0\circ\Lambda_0>\theta\).

It remains to verify the condition (iii) in part (b). Since
$\alpha_0=\hat\alpha+\hat\gamma$ on
$\what\Theta_2\sm\hat\Lambda_0\inv(\what V_2\cap\what
Y)\subset\what\Theta_2\sm\what Q$, since $\alpha+\gamma$ is of
class $\strplshclass^\infty(g,q)$ on a \nbd of
$B\cap\overline\Theta_2\sm Q$, and since $\gamma$ is \cinf \str
\plsh on a \nbd of $\overline V_1\cap\overline\Theta_2$, we need
only show that $\alpha_1$ is \cinf \strg $q$-\cvx in a \nbd of
each point
$x\in\what\Theta_2\cap\hat\Lambda_0\inv\left(\overline{\what
V_2}\cap\what Y\right)$.  For this, observe that, since
\(y=\hat\Lambda_0(x)\in\what\Omega_1\cap\what Y\), we have
$\lap_{\hat g}(\hat\beta+\beta_0)(y)>0$. Hence there exists a
$1$-dimensional vector subspace $\cal V$ of $T^{1,0}_y\what Y$ \st
$\lev{\hat\beta+\beta_0}>0$ on $\cal V\sm\set{0}$. Since
$\Lambda_0$ is a submersion, the inverse image $\cal
U=[(\Lambda_0)_*]_x\inv (\cal V)$ is a vector subspace of
dimension $n-q+1$ in $T^{1,0}_y\what X$. By construction, on some
\nbd of $x$ we have
$\hat\beta\circ\hat\Lambda_0+\hat\lambda\hat\rho>3\delta$ and
hence
\[
\alpha_1=\chi(\hat\beta\circ\hat\Lambda_0+\hat\lambda\hat\rho)+\beta_0\circ\hat\Lambda_0
=\hat\beta\circ\hat\Lambda_0+\hat\lambda\hat\rho-2\delta+\beta_0\circ\hat\Lambda_0
=(\hat\beta+\beta_0)\circ\hat\Lambda_0+\hat\lambda\hat\rho-2\delta.
\]
Moreover, $\hat\lambda$ is constant near $x$. Thus for each
tangent vector $v\in\cal U\sm\set{0}$, we have
\[
\lev{\alpha_1}(v,v)=\lev{\hat\beta+\beta_0}\left((\hat\Lambda_0)_*v,(\hat\Lambda_0)_*v\right)+\hat\lambda(x)\lev{\hat\rho}(v,v).
\]
Both terms on the right-hand side are nonnegative, the first is
positive if $(\hat\Lambda_0)_*v\neq 0$, and the second is positive
if $(\hat\Lambda_0)_*v=0$. Thus $\alpha_1$ is \cinf \strg $q$-\cvx
near $x$.

We now consider the case of a general \con reduced \cpx space $X$
of bounded local \holo embedding dimension. According to the Stein
\nbd theorem of Siu \cite{Siu-Stein neighborhoods} and the
embedding theorem of Remmert~\cite{Remmert embedding},
Narasimhan~\cite{Narasimhan embedding}, and Bishop~\cite{Bishop
Stein embedding}, some \con \nbd of $Y$ in $X$ admits a proper
\holo embedding into $\C^n$ for some positive integer $n$. Thus
there exists a proper \holo embedding $\Phi\colon Z\to X'$ of a
\con \nbd $Z$ of $Y$ in $X$ into a \con open subset $X'$ of $\C^n$
\st $Y'=\Phi(Y)$ is closed in $\C^n$, $Y$ is a \cont strong
deformation retract of $Z$, and $Z'=\Phi(Z)$ is a \cont strong
deformation retract of $X'$. We may also fix a Hermitian metric
$g'$ on $X'$ \st $\Phi^*g'=g$ on $Z$. Finally, we may fix a \con
\nbd $P_1$ of $Y'$ in $\C^n$ with $\overline P_1\subset X'$.

Given open sets $\Omega_1$, $\Omega_2$, and $Q$ as in the
statement of the proposition, we may choose open sets $\Omega'_1$,
$\Omega'_2$, and $Q'$ in $\C^n$ \st
$X'\Supset\Omega'_1\Supset\Omega'_2\Supset Q'$, $Q'\cap Z'\cap
P_1=\Phi(Q\cap Z)\cap P_1$, $\Omega'_j\cap Z'\cap
P_1=\Phi(\Omega_j\cap Z)\cap P_1$ for $j=1,2$, $\overline{Q'\cap
Z'}=\overline{Q'}\cap Z'$, and $\overline{\Omega'_j\cap
Z'}=\overline{\Omega'_j}\cap Z'$ for $j=1,2$. By the above, there
exists a \nbd $\Theta_1'\subset P_1$ of $Y'$ and a \cinf \fn
$\alpha'$ satisfying the conditions in part~(a) relative to these
choices of sets in $\C^n$. The \concomp $\Theta_1$ of
$\Phi\inv(\Theta'_1)$ containing $Y$ and the \fn
$\alpha=\alpha'\circ\Phi\restrict{\Theta_1}$ on $\Theta_1$ then
satisfy the conditions in part~(a) in $X$.

Fix a \con open set $P_2$ in $\C^n$ with $\overline
P_2\subset\Theta_1$. Given open sets $\Omega_3$ and $\Omega_4$ as
in part~(b), we may choose open sets $\Omega'_3$ and $\Omega'_4$
in $\C^n$ \st $\Omega'_2\Supset\Omega'_3\Supset\Omega'_4\Supset
Q'$, $\Omega'_j\cap Z'\cap P_2=\Phi(\Omega_j\cap Z)\cap P_2$ for
$j=3,4$, and $\overline{\Omega'_j\cap Z'}=\overline{\Omega'_j}\cap
Z'$ for $j=3,4$. Set $B'=\setof{x\in\Theta'_1}{\alpha'(x)>0}$ and
$B=\Phi\inv(B')\cap\Theta_1=\setof{x\in\Theta_1}{\alpha(x)>0}$.
Given a family of \fns $\seq\omega\nu_{\nu\in N}$ as in the
statement of part~(b), we may form a family of real-valued \cinf
\fns $\set{\omega'_\nu}_{\nu\in N}$ with \cpt support in $B'$ \st
$\omega'_\nu\circ\Phi=\omega_\nu$ on $B$ for each~$\nu$.

Again, for some \nbd $\Theta_2'\subset P_2$ of $Y'$ and for
sufficiently large positive constants $\seq R\nu$, we get the
conditions in part~(b) relative to the above choices and any
covering space of $X'$ (with the required properties). Let
$\Theta_2$ be the \concomp of $\Phi\inv(\Theta'_2)$ containing
$Y$. For any infinite covering space $\Upsilon\colon\what X\to X$
in which $\what Y=\Upsilon\inv (Y)$ and $\Upsilon\inv
(\Omega_4\cap Y)$ are \con and for any positive \cont \fn $\theta$
on $\what X$, we get a corresponding covering
$\Upsilon'\colon\what X'\to X'$ \st $\Phi$ lifts to a proper \holo
embedding $\hat\Phi$ of $\what Z=\Upsilon\inv (Z)$ into $\what X'$
and we may choose a positive \cont \fn $\theta'$ on $\what X'$
with $\theta'\circ\hat\Phi=\theta$ on $\what Z$. Forming the
corresponding \fn $\alpha'_0$ on
$\what\Theta'_2=(\Upsilon')\inv(\Theta_2')$, we get the desired
\fn $\alpha_0=\alpha'_0\circ\hat\Phi\restrict{\what\Theta_2}$ on
$\what\Theta_2=\Upsilon\inv(\Theta_2)$.
\end{pf*}

\section{A uniformly quasi-plurisubharmonic function}\label{Quasi-plsh fn near singular set section}

We recall that an upper semi-\cont \fn $\vphi\colon
M\to[-\infty,\infty)$ on a \cpx manifold $M$ is called
\textit{\plshns} if, for every \holo mapping $f\colon D\to M$ of a
disk $D\subset\C$ into $M$, the \fn $\psi=\vphi\circ f$ is
subharmonic in~$D$; that is, for every \cont \fn $u$ on a \cpt set
$K\subset D$ which is \harm on the interior $\overset\circ K$ and
which satisfies $u\geq\psi$ on $\partial K$, we have $u\geq\psi$
on $K$. The \fn $\vphi$ is \textit{\str \plshns} if, for every
real-valued \cinf \fn $\rho$ with \cpt support in $M$, the \fn
$\vphi+\epsilon\rho$ is \plsh for every sufficiently small
$\epsilon>0$. Let $k\in\Z_{\geq 0}\cup\set{\infty,\omega}$. A \fn
$\vphi\colon X\to[-\infty,\infty)$ on a reduced \cpx space~$X$ is
\textit{\plshns} (\textit{\str \plshns, $C^k$ \plshns, $C^k$ \str
\plshns}) if, for each point in $p\in X$, there is a local \holo
chart $(U,\Phi,U')$ with $p\in U$ and a \fn $\vphi'$ on $U'$ \st
$\vphi'$ is \plsh (respectively, \str \plshns, \plsh and of class
$C^k$, \str \plsh and of class $C^k$) and \st
$\vphi=\vphi'\circ\Phi$ on $U$. The \fn $\vphi$ is
\textit{quasi-\plshns} (\textit{$C^k$ quasi-\plshns}) if, locally,
$\vphi$ is the sum of a real \anal \fn and a \plsh (respectively,
$C^k$ \plshns) \fnns.  Clearly, for $k\in\Z_{\geq
2}\cup\set{\infty,\omega}$, a $C^k$ \fn is $C^k$ quasi-\plshns. By
a theorem of Richberg~\cite{Richberg}, a \cont \fn which is \str
\plsh is automatically $C^0$ \str \plshns. However, the
corresponding property for $C^k$ \fns does not always hold on a
\cpx space (see, for example, Smith~\cite{Smith-cinf plsh ctex}).
Finally, we observe that the class of \cinf strongly $1$-\cvx \fns
is precisely the class of \cinf \str \plsh \fns and
$\plshclass^\infty(1)$ is precisely the class of $\cinf$ \plsh
\fns (see the introduction).

The goal of this section is to produce a \fn on each covering
space which is, in a sense, uniformly quasi-\plsh and which has a
logarithmic singularity along the lifting of a given \cpt \anal
set (cf.~Lemma~5 of Demailly~\cite{Demailly-Cohomology of q-convex
spaces}). For this, it will be convenient to have the following
terminology:
\begin{defn}\label{Quasi-plsh with sing def}
Let $A$ be a (properly embedded) \anal subset of a reduced \cpx
space~$X$ and let $\vphi\colon\Omega\to[-\infty,\infty)$ be a
mapping on an open subset $\Omega$ of $X$. We will say that
$\vphi$ is of class $\plshclass^\infty_A$
($\strplshclass^\infty_A$, $\cal Q^\infty_A$) if, for each point
$p\in\Omega$, there is a local \holo model $(U,\Phi,U')$ in $X$
with $p\in U\subset\Omega$ and a \plsh (respectively, \str
\plshns, quasi-\plshns) \fn $\vphi'\colon U'\to[-\infty,\infty)$
\st $\vphi=\vphi'\circ\Phi$ on $U$, $\vphi'$ is real-valued and
\cinf on $U'\sm\Phi(A\cap U)$, $\vphi'$ is \cont on $U'$ as a
mapping into $[-\infty,\infty)$, and $\vphi'=-\infty$ on
$\Phi(A\cap U)$.
\end{defn}

The main goal of this section is the following:
\begin{prop}\label{Quasi-plsh fn near the singular set prop}
Let $(X,g)$ be a \con reduced Hermitian \cpx space, let $Y$ be a
\cpt \anal subset of~$X$, let $\Omega$ be a \rel \cpt \nbd of~$Y$
in~$X$, let $\seq D\nu_{\nu=1}^l$ be \rel \cpt open subsets
of~$X$, and let $\rho_\nu$ be a \cinf \str \plsh \fn on a \nbd of
$\overline D_\nu$ for each $\nu=1,\dots,l$. Then, for some choice
of nondecreasing \cont maps
\[
\eta_-\colon[0,\infty]\to[-\infty,0]\qquad\text{and}\qquad\eta_+\colon
[0,\infty]\to [-\infty,0]
\]
with \(\eta_-\leq\eta_+\), \(\eta_+(0)=-\infty\), and
\(\eta_->-\infty\) on $(0,\infty]$; for every sufficiently large
$R>0$; for every \con covering space $\Upsilon\colon\what X\to X$;
and for every \anal set $Z$ equal to a union of \ircomps of $\what
Y\equiv\Upsilon\inv(Y)$; there exists a \cont \fn
$\alpha\colon\what X\to [-\infty,0]$ with the following
properties:
\begin{enumerate}
\item[(i)] The \fn $\alpha$ is of class $\cal Q^\infty_Z$ on
$\what X$;

\item[(ii)] We have
$\supp\alpha\subset\what\Omega\equiv\Upsilon\inv(\Omega)$;

\item[(iii)] For each $\nu=1,\dots,l$,
$\alpha+R\cdot\rho_\nu\circ\Upsilon$ is of class
$\strplshclass^\infty_Z$ on $\Upsilon\inv(D_\nu)$;

\item[(iv)] For $\hat g=\Upsilon^*g$, we have
\[
\eta_-(\dist_{\hat g}(x,Z))\leq\alpha(x)\leq\eta_+(\dist_{\hat
g}(x,Z))\qquad\forall\, x\in\what X.
\]

\end{enumerate}
\end{prop}
\begin{rmks} 1. For $Z=\emptyset$, we have $\dist_{\hat
g}(\cdot,Z)\equiv\infty$.

\noindent 2. Clearly, we may choose $R=1$ (by replacing $\alpha$
and $\eta_{\pm}$ with $R\inv\alpha$ and $R\inv\eta_{\pm}$,
respectively), but it will be more convenient for the proof to
allow $R$ to vary.

\noindent 3. The condition $\strplshclass^\infty_Z$ is slightly
stronger than necessary for our purposes, but working with this
class allows one to write some of the arguments more efficiently.
\end{rmks}

The proof of Proposition~\ref{Quasi-plsh fn near the singular set
prop} is a modification of Demailly's construction of a
quasi-\plsh \fn with logarithmic singularities along a given \anal
subset (Lemma~5 of \cite{Demailly-Cohomology of q-convex spaces}).
The \fn is obtained by patching local quasi-\plsh \fns using the
following \cinf version of the maximum \fn for which the given
properties are easy to check:

\begin{lem}\label{cinf max lemma}
Let $\kappa\colon\R\to[0,\infty)$ be a $\cinf$ \fn \st
$\supp\kappa\subset (-1,1)$, $\int_{\R}\kappa (u)\,du=1$, and
$\int_{\R}u\kappa (u)\,du=0$. For each $d\in\pos{\Z}$, let $\cal
M_d\colon\R^d\to\R$ be the \fn given by
\[
\cal M_d(t)=\int_{\R^d}\biggl[\max_{1\leq j\leq
d}(t_j+u_j)\biggr]\prod_{1\leq j\leq d}\kappa (u_j)\,du_j
\]
for each $t=(t_1,\dots,t_d)\in\R^d$. Then
\begin{enumerate}

\item[(a)] For each $t=(t_1,\dots,t_d)\in\R^d$,
\[
\cal M_d(t)=\int_{\R^m}\biggl[\max_{1\leq j\leq
d}u_j\biggr]\prod_{1\leq j\leq d}\kappa(u_j-t_j)\,du_j.
\]

\item[(b)] $\cal M_d$ is $\cinf$, \cvxns, and symmetric in
$t_1,\dots,t_d$.

\item[(c)] $\cal M_d(t_1,\dots,t_d)$ is nondecreasing in each
variable $t_j$.

\item[(d)] For every $t=(t_1,\dots,t_d)\in\R^d$ and $s\in\R$, we
have
\[
\cal M_d(t_1+s,\dots,t_d+s)=\cal M_d(t)+s.
\]

\item[(e)] For every $t\in\R^d$, \(\max(t)\leq\cal
M_d(t)\leq\max(t)+1\).

\item[(f)] If $t'=(t_0,t_1,\dots,t_d)=(t_0,t)\in\R^{d+1}$ with
$t_0\leq t_1-2$, then $\cal M_{d+1}(t')=\cal M_d(t)$.

\item[(g)] For $d=1$, $\cal M_d(t)=t$ for each $t\in\R$.

\item[(h)] For each $j=1,\dots,d$, we have
$0\leq\frac{\partial}{\partial t_j}\cal M_d(t_1,\dots,t_d)\leq 1$.

\item[(i)] If $\vphi=(\vphi_1,\dots,\vphi_d)$ is a $d$-tuple of
\cinf real-valued \fns on a \cpx manifold~$X$, then
\[
\lev{\cal M_d(\vphi)}(v,v)\geq\min_{1\leq j\leq
d}\lev{\vphi_j}(v,v)\qquad\forall\, v\in T^{1,0}X.
\]
\end{enumerate}
\end{lem}

The following is Demailly's construction:
\begin{lem}[See Lemma~5 of \cite{Demailly-Cohomology of q-convex spaces}]\label{quasiplsh fn on the base lemma}
Let $X$ be a reduced \cpx space, let $A$ be an \anal subset of
$X$, let $\seq Uj_{j\in J}$ and $\seq Vj_{j\in J}$ be locally
finite coverings of $X$ by open sets with $\overline V_j\subset
U_j$ for each $j\in J$, and, for each index $j\in J$, let
$\lambda_j$ be a \cinf \fn on~$V_j$ with $\lambda_j\to-\infty$ at
$\partial V_j$ and let
\[
F_j=\left(f_j\ssp 1,\dots,f_j\ssp{N_j}\right)\colon U_j\to\C^{N_j}
\]
be a \holo map \st the coordinate \fns generate the ideal sheaf
$\ideal {(A\cap U_j)}$ at each point in $U_j$. Let $\alpha\colon
X\to[-\infty,\infty)$ be the mapping defined by setting
$\alpha\equiv-\infty$ on~$A$ and, for $x\in X\sm A$, setting
\[
\alpha(x)=\cal
M_d\left[\left(\log\left|F_{j_1}(x)\right|^2+\lambda_{j_1}(x)\right),\dots,
\left(\log\left|F_{j_d}(x)\right|^2+\lambda_{j_d}(x)\right)\right];
\]
where $j_1,\dots,j_d$ are the distinct indices with $x\in V_j$ if
and only if $j\in\set{j_1,\dots,j_d}$. Then $\alpha$ is of class
$\cal Q^\infty_A$ on $X$.
\end{lem}
\begin{pf}
Given a point $x\in X$, we have distinct indices
$j_1,\dots,j_d,j_{d+1},\dots,j_m\in J$ \st $x\in V_j$
(respectively, $x\in\overline V_j$) if and only if
$j\in\set{j_1,\dots,j_d}$ (respectively,
$j\in\set{j_1,\dots,j_m}$). Thus we may choose a \nbd $W$ of $x$
with
\[
W\Subset V_{j_1}\cap\cdots\cap V_{j_d}\cap
U_{j_{d+1}}\cap\cdots\cap U_{j_m}\text{ and }W\cap
V_j=\emptyset\quad\forall\, j\in J\sm\set{j_1,\dots,j_m}.
\]
Clearly, the \fns $|F_{j_\mu}|/|F_{j_\nu}|$ are bounded on
$\overline W$ for $\mu,\nu=1,\dots,m$. On the other hand,
$\lambda_{j_\nu}$ is bounded on $\overline W$ for $\nu=1,\dots,d$
and $\lambda_{j_\nu}\to-\infty$ at $x\in\partial V_{j_\nu}$ for
$\nu=d+1,\dots,m$. Hence, choosing $W$ sufficiently small, we get
\[
\log\left|F_{j_\mu}\right|^2+\lambda_{j_\mu}-2>\log\left|F_{j_\nu}(x)\right|^2+\lambda_{j_\nu}
\]
on $W\cap V_{j_\nu}$ for $1\leq\mu\leq d$ and $d+1\leq\nu\leq m$.
Thus
\[
\alpha=\cal
M_d\left[\left(\log\left|F_{j_1}\right|^2+\lambda_{j_1}\right),\dots,
\left(\log\left|F_{j_d}\right|^2+\lambda_{j_d}\right)\right]
\]
on $W\sm Z$.

Now we may also choose $W$ so that we have a proper local \holo
model $(W,\Phi,W')$ and, for each $\nu=1,\dots,d$, a \holo mapping
$F'_{j_\nu}\colon W'\to\C^{N_{j_\nu}+l_{j_\nu}}$ \st
$F'_{j_\nu}\circ\Phi=(F_{j_\nu},0,\dots,0)$ on $W$ and \st the
coordinate \fns for $F'_{j_\nu}$ generate the ideal sheaf
$\ideal{\Phi(A\cap W)}$ at each point. Furthermore, we may assume
that there is a \cinf \str \plsh \fn $\rho$ on $W'$ and \cinf \fns
$\lambda'_{j_1},\dots,\lambda'_{j_d}$ on $W'$ \stns, for each
$\nu=1,\dots,d$, we have
$\lambda'_{j_\nu}\circ\Phi=\lambda_{j_\nu}$ on $W$ and
$\lambda'_{j_\nu}+\rho$ is \plsh on $W'$.  Setting
$\alpha'\equiv-\infty$ on $\Phi(A\cap W)$ and
\[
\alpha'=\cal
M_d\left[\left(\log\left|F'_{j_1}\right|^2+\lambda'_{j_1}\right),\dots,
\left(\log\left|F'_{j_d}\right|^2+\lambda'_{j_d}\right)\right]
\]
on $W'\sm\Phi(A\cap W)$, we see that $\alpha'\circ\Phi=\alpha$ on
$W$, $\alpha'$ is \cinf on $W'\sm\Phi(A\cap W)$, and $\alpha'$ is
\cont as a mapping into $[-\infty,\infty)$. Moreover, by part~(i)
of Lemma~\ref{cinf max lemma}, $\alpha'+\rho$ is \plsh on
$W'\sm\Phi(A\cap W)$ and, therefore, on $W'$. It follows that
$\alpha$ is of class $\cal Q^\infty_A$ on $X$.
\end{pf}

\begin{pf*}{Proof of Proposition~\ref{Quasi-plsh fn near the singular set
prop}} We may fix finite coverings $\seq Qi_{i=1}^k$, $\seq
Ui_{i=1}^k$, and $\seq Vi_{i=1}^k$ of $Y$ by nonempty \con open
subsets of $X$ and \fns $\seq\lambda i_{i=1}^k$ \stns, for each
$i=1,\dots,k$,
\begin{enumerate}
\item[(\ref{Quasi-plsh fn near the singular set prop}.1)] $Q_i$ is
contractible;

\item[(\ref{Quasi-plsh fn near the singular set prop}.2)]
$V_i\Subset U_i\Subset Q_i\Subset\Omega$;

\item[(\ref{Quasi-plsh fn near the singular set prop}.3)] We have
$\diam_gU_i<\dist_g(U_i,X\sm Q_i)$ (where the diameter $\diam_g$
is taken with respect to the distance \fn $\dist_g(\cdot,\cdot)$
in $X$;

\item[(\ref{Quasi-plsh fn near the singular set prop}.4)] If
$\overline{U_i}\cap \overline{U_j}\neq\emptyset$ for some~$j$,
then $U_i\Subset Q_j$;

\item[(\ref{Quasi-plsh fn near the singular set prop}.5)] For
every set $A$ which is equal to a union of \ircomps of $Y\cap Q_i$
each of which meets $U_i$, we have a \holo map
\[
F_{i,A}=\left(f_{i,A}\ssp
1,\dots,f_{i,A}\ssp{N_{i,A}}\right)\colon U_i\to\C^{N_{i,A}}
\]
\st the coordinate \fns generate the ideal sheaf $\ideal A$ at
each point in~$U_i$ and
\[
|F_{i,A}|^2=\sum_{j=1}^{N_{i,A}}|f_{i,A}\ssp j|^2\leq 1
\]
(here, we allow for $A=\emptyset$ in which case we take
$N_{i,A}=1$ and $F_{i,A}=f_{i,A}\ssp 1\equiv 1$); and

\item[(\ref{Quasi-plsh fn near the singular set prop}.6)] We have
$\lambda_i\in\cinf(V_i)$, $\lambda_i<-2$ on~$V_i$, and
$\lambda_i\to-\infty$ at $\partial V_i$.

\end{enumerate}
We may also choose open sets $Q_0$, $U_0$, $V_0$, and $W_0$ in $X$
with
\[
X\sm (V_1\cup\cdots\cup V_k)\subset
W_0\subset\overline{W_0}\subset V_0\subset\overline{V_0}\subset
U_0\subset\overline{U_0}\subset Q_0\subset\overline{Q_0}\subset
X\sm Y
\]
and a nonpositive \fn $\lambda_0\in\cinf (V_0)$ with
$\lambda_0\equiv 0$ on $W_0$ and $\lambda_0\to-\infty$
at~$\partial V_0$. We also set $N_{0,\emptyset}=1$ and
\[
F_{0,\emptyset}=f_{0,\emptyset}\ssp
1=f_{0,\emptyset}\ssp{N_{0,\emptyset}}\equiv 1.
\]

For $i=0,1,2,\dots,k$, let $\cal A_i$ be the (finite) collection
of \anal subsets of $Q_i$ which are either empty or the union of a
set of \ircomps of $Y\cap Q_i$ each of which meets~$U_i$. Given an
open subset $\Theta$ of $X$, let $\cal A_\Theta$ be the collection
of $(k+1)$-tuples $A=(A_0,\dots,A_k)\in\cal
A_0\times\cdots\cal\times A_k$ \stns, for all $i,j=0,\dots,k$, we
have
\[
A_i\cap U_i\cap U_j\cap\Theta=A_j\cap U_j\cap U_i\cap\Theta.
\]
In particular, the associated set $\check A\equiv[(A_0\cap
U_0)\cup\cdots\cup (A_k\cap U_k)]\cap\Theta$ is a (properly
embedded) \anal subset of $\Theta$ with $\check A\cap U_i=A_i\cap
U_i\cap\Theta$ for $i=0,\dots,k$. Thus for each
$A=(A_0,\dots,A_k)\in\cal A_\Theta$, as in Lemma~\ref{quasiplsh fn
on the base lemma}, we may define a \cont \fn
$\alpha_{\Theta,A}\colon\Theta\to[-\infty,0]$ of class $\cal
Q^\infty_{\check A}$ by setting $\alpha_{\Theta,A}\equiv-\infty$
on $\check A$ and, for each $x\in\Theta\sm\check A$, setting
\[
\alpha_{\Theta,A}(x)=\cal
M_d\left[\left(\log\left|F_{i_1,A_{i_1}}(x)\right|^2+\lambda_{i_1}(x)\right),\dots,
\left(\log\left|F_{i_d,A_{i_d}}(x)\right|^2+\lambda_{i_d}(x)\right)\right];
\]
where $i_1,\dots,i_d$ are the distinct indices with $x\in V_i$ if
and only if $i\in\set{i_1,\dots,i_d}$. We have
$\alpha_{\Theta,A}\leq 0$ because, in the above notation, all of
the entries are nonpositive and at most one is greater than~$-2$;
while, if $t=(t_1,\dots,t_d)$ with $t_1\leq 0$ and $t_j<-2$ for
$j=2,\dots,d$, then
\[
\cal M_d(t)\leq\cal M_d(0,-2,\dots,-2)=\cal M_1(0)=0.
\]
Similarly, we have $\alpha_{\Theta,A}\equiv 0$ on
$W_0\cap\Theta\supset\Theta\sm\Omega$.  Since $X\sm
W_0\subset\Omega\Subset X$, it follows that, given an open set $V$
with $\overline V\subset\Theta$, there exists a positive
constant~$R_{\Theta,A,V}$ and nondecreasing \cont \fns
\[
\eta_{\Theta,A,V,-}\colon
[0,\infty]\to[-\infty,0]\qquad\text{and}\qquad\eta_{\Theta,A,V,+}\colon
[0,\infty]\to[-\infty,0]
\]
\st \(\eta_{\Theta,A,V,-}\leq\eta_{\Theta,A,V,+}\),
\(\eta_{\Theta,A,V,+}(0)=-\infty\),
\(\eta_{\Theta,A,V,-}>-\infty\) on $(0,\infty]$,
\[
\eta_{\Theta,A,V,-}(\dist_{g}(x,\check A))\leq\alpha_{\Theta,A}(x)
\leq\eta_{\Theta,A,V,+}(\dist_{g}(x,\check A))\qquad\forall\, x\in
V,
\]
and, for each $\nu=1,\dots,l$, the \fn
$\alpha_{\Theta,A}+R_{\Theta,A,V}\cdot\rho_\nu$ is of class
$\strplshclass^\infty_{\check A}$ on $D_\nu\cap V$.

The collection $\cal A_{U_i}\subset\cal A_0\times\cdots\times\cal
A_k$ is finite for each~$i$. Thus we may choose a constant $R_0>0$
and nondecreasing \cont \fns
\[
\eta_-\colon
[0,\infty]\to[-\infty,0]\qquad\text{and}\qquad\eta_+\colon
[0,\infty]\to[-\infty,0]
\]
\st \(\eta_-\leq\eta_+\), \(\eta_+(0)=-\infty\),
\(\eta_->-\infty\) on $(0,\infty]$, and, for each $i=0,1,\dots,k$
and each $A\in\cal A_{U_i}$, we have $R_0>R_{U_i,A,V_i}$,
\(\eta_-\leq\eta_{U_i,A,V_i,-}\leq\eta_{U_i,A,V_i,+}\leq\eta_+\),
and $\eta_+\equiv 0$ on the interval $[\dist_g(V_i,X\sm
U_i),\infty]$ (note that $\dist_g(V_0,X\sm U_0)>0$ since $X\sm
U_0$ is \cptns).

Given a \con covering space $\Upsilon\colon\what X\to X$ and an
\anal subset $Z$ of $X$ which is equal to a union of \ircomps of
$\what Y=\Upsilon\inv(Y)$, we will construct the associated class
$\cal Q^\infty_Z$ \fn $\alpha\colon\what X\to[-\infty,0]$. For
this, we assume that the covering is infinite.  The proof for a
finite covering is similar.

For each $i=0,\dots,k$, we let $\what Q_i=\Upsilon\inv(Q_i)$,
$\what U_i=\Upsilon\inv(U_i)$, and $\what V_i=\Upsilon\inv(V_i)$.
We set $\what W_0=\Upsilon\inv(W_0)$. For each $i=1,\dots,k$, we
let $\set{Q\ssp\nu_i}_{\nu=1}^\infty$,
$\set{U\ssp\nu_i}_{\nu=1}^\infty$, and
$\set{V\ssp\nu_i}_{\nu=1}^\infty$ denote the distinct \concomps of
$\what Q_i$, $\what U_i$, and $\what V_i$, respectively, with
$V\ssp\nu_i\Subset U\ssp\nu_i\Subset Q\ssp\nu_i$ for each~$\nu$.
For $1\leq i\leq k$ and $\nu\geq 1$, we let
$\lambda\ssp\nu_i=\lambda_i\circ\Upsilon\restrict{V\ssp\nu_i}\in\cinf\left(V\ssp\nu_i\right)$,
we let $Z\ssp\nu_i$ be the union of all \ircomps of~$Z\cap
Q\ssp\nu_i$ meeting $U\ssp\nu_i$, we let
$A\ssp\nu_i=\Upsilon\left(Z\ssp\nu_i\right)\in\cal A_i$, and we
let
$F\ssp\nu_i=F_{i,A\ssp\nu_i}\circ\Upsilon\restrict{U\ssp\nu_i}\colon
U\ssp\nu_i\to\C^{N_{i,A\ssp\nu_i}}$. We also set $Q\ssp 1_0=\what
Q_0$, $U\ssp 1_0=\what U_0$, $V\ssp 1_0=\what V_0$, $W\ssp
1_0=\what W_0$, $\lambda\ssp
1_0=\lambda_0\circ\Upsilon\in\cinf\left(V\ssp 1_0\right)$, $Z\ssp
1_0=\emptyset$, $A\ssp 1_0=\emptyset$, and $F\ssp
1_0=F_{0,\emptyset}\circ\Upsilon\equiv 1$ on $U\ssp 1_0$.

Lemma~\ref{quasiplsh fn on the base lemma} applied to the above
objects yields a \fn $\alpha\colon X\to[-\infty,0]$ satisfying the
conditions~(i)~and~(ii). Specifically, $\alpha\equiv-\infty$ on
$Z$ while, for $x\in\what X\sm Z$ and for distinct pairs of
indices $(i_1,\nu_1),\dots,(i_d,\nu_d)$ with $x\in V\ssp\nu_i$ if
and only if $(i,\nu)\in\set{(i_1,\nu_1),\dots,(i_d,\nu_d)}$, we
have,
\[
\alpha(x)=\cal
M_d\left[\left(\log\left|F\ssp{\nu_1}_{i_1}(x)\right|^2+\lambda\ssp{\nu_1}_{i_1}(x)\right),\dots,
\left(\log\left|F\ssp{\nu_d}_{i_d}(x)\right|^2+\lambda\ssp{\nu_d}_{i_d}(x)\right)\right].
\]
We will show that $\alpha$ satisfies the conditions~(iii) (for
$R>R_0$) and (iv) by showing that, locally, we have
$\alpha=\alpha_{U_i,A}\circ\Upsilon$ with $A\in\cal A_{U_i}$.

Given $i\in\set{1,\dots,k}$ and $\nu\in\N$, we have
\textit{distinct} indices $i_1,\dots,i_m\in\set{0,\dots,k}$ and
indices $\nu_1,\dots,\nu_m\in\N$ \stns, for any pair of indices
$(j,\mu)$, we have
\[
U\ssp{\nu}_{i}\cap
U\ssp\mu_j\neq\emptyset\iff(j,\mu)\in\set{(i_1,\nu_1),\dots,(i_m,\nu_m)}
\]
(here we have used the condition (\ref{Quasi-plsh fn near the
singular set prop}.4)). In particular,
$(i,\nu)\in\set{(i_1,\nu_1),\dots,(i_m,\nu_m)}$. For each
$j=0,\dots,k$, let $A_j\in\cal A_j$ be the \anal set given by
\[
A_j=\left\{
\begin{aligned}
\emptyset&\text{ if }j\not\in\set{i_1,\dots,i_m}\\
A\ssp{\nu_s}_{i_s}&\text{ if }j=i_s
\end{aligned}
\right.
\]
Then the $k$-tuple $A\equiv (A_0,\dots,A_k)$ is an element of
$\cal A_{U_i}$ with associated \anal set
\[
\check A=\Upsilon\left(Z\cap
U\ssp\nu_i\right)=\Upsilon\left(Z\ssp\nu_i\cap
U\ssp\nu_i\right)=A\ssp\nu_i\cap U_i=A_i\cap U_i.
\]
To see this, suppose $0\leq j\leq k$. If $j=0$ or
$j\notin{i_1,\dots,i_m}$, then \( A_j\cap U_j\cap
U_i=\emptyset=A_i\cap U_i\cap U_j\). If $j=i_s\neq 0$ for some
$s$, then we have
\[
A_j\cap U_j\cap U_i=A\ssp{\nu_s}_{i_s}\cap U_{i_s}\cap
U_i=\Upsilon\left(Z\ssp{\nu_s}_{i_s}\cap
U\ssp{\nu_s}_{i_s}\right)\cap U_i =\Upsilon\left(Z\cap
U\ssp{\nu_s}_{i_s}\cap U\ssp\nu_i\right);
\]
where the last equality holds because $U\ssp\nu_i$ is the unique
\comp of $\what U_i$ meeting $U\ssp{\nu_s}_{i_s}$. Exchanging $j$
and $i$, we see that \( A_j\cap U_j\cap U_i=A_i\cap U_i\cap U_j \)
in this case as well and the claim follows.

Furthermore, we have $\alpha=\alpha_{U_i,A}\circ\Upsilon$ on
$U\ssp\nu_i$. For both \fns equal $-\infty$ on $Z\cap U\ssp\nu_i$.
Given a point $x\in U\ssp\nu_i\sm Z$, we may rearrange the indices
so that, for some $d\in{1,\dots,m}$, we have $x\in
V\ssp{\nu_1}_{i_1}\cap\cdots\cap V\ssp{\nu_d}_{i_d}$ and
$x\notin\what V_j$ for $j\notin\set{i_1,\dots,i_d}$. Hence
\[
\alpha(x)=\cal
M_d\left[\left(\log\left|F\ssp{\nu_1}_{i_1}(x)\right|^2+\lambda\ssp{\nu_1}_{i_1}(x)\right),\dots,
\left(\log\left|F\ssp{\nu_d}_{i_d}(x)\right|^2+\lambda\ssp{\nu_d}_{i_d}(x)\right)\right].
\]
We also have $y=\Upsilon(x)\in V_{i_1}\cap\cdots\cap V_{i_d}$ and
$y\notin V_j$ for $j\notin\set{i_1,\dots,i_d}$. Thus
\[
\alpha_{U_i,A}(y)=\cal
M_d\left[\left(\log\left|F_{i_1,A_{i_1}}(y)\right|^2+\lambda_{i_1}(y)\right),\dots,
\left(\log\left|F_{i_d,A_{i_d}}(y)\right|^2+\lambda_{i_d}(y)\right)\right];
\]
and we get the claim.

It follows that the condition~(iii) holds on the set $V\ssp\nu_i$
for $R>R_0$. That is, for $R>R_0$, the \fn
$\alpha+R\cdot\rho_\nu\circ\Upsilon$ is of class
$\strplshclass^\infty_Z$ on $V\ssp\nu_i\cap\Upsilon\inv(D_\mu)$
for $\mu=1,\dots,l$. To verify that the condition~(iv) also holds
on $V\ssp\nu_i$, we fix $x\in V\ssp\nu_i$ and set $y=\Upsilon(x)$
and $r=\dist_{\hat g}(x,Z)$, where $\hat g=\Upsilon^*g$. The
condition~(\ref{Quasi-plsh fn near the singular set prop}.3)
ensures that $\dist_{\hat
g}(a,b)=\dist_g(\Upsilon(a),\Upsilon(b))$ for all $a,b\in
U\ssp\nu_i$ and
\[
\diam_{\hat g}U\ssp\nu_i=\diam_gU_i<\dist_g(U_i,X\sm
Q_i)=\dist_{\hat g}\left(U\ssp\nu_i,\what X\sm Q\ssp\nu_i\right).
\]
In particular, we have
\begin{align*}
\eta_-(r)&\leq\eta_-\left(\dist_{\hat g}\left(x,Z\cap
U\ssp\nu_i\right)\right)=\eta_-(\dist_g(y,\check
A))\\
&\leq\eta_{U_i,A,V_i,-}(\dist_g(y,\check A))
\leq\alpha_{U_i,A}(y)=\alpha(x)
\end{align*}
For the other required inequality, we observe that, if
\[
r\geq\dist_{\hat g}\left(V\ssp\nu_i,\what X\sm
U\ssp\nu_i\right)=\dist_g(V_i,X\sm U_i),
\]
then, by the choice of
$\eta_+$, we have
\[
\alpha(x)=\alpha_{U_i,A}(y)\leq 0=\eta_+(r).
\]
If $r<\dist_{\hat g}(V\ssp\nu_i,\what X\sm U\ssp\nu_i)$, then we
must have
\[
r=\dist_{\hat g}\left(x,Z\cap U\ssp\nu_i\right)=\dist_g(y,\check
A)
\]
and hence
\[
\alpha(x)=\alpha_{U_i,A}(y)\leq\eta_{U_i,A,V_i,+}(r)\leq\eta_+(r).
\]
Thus the condition~(iv) holds on $V\ssp\nu_i$.

It remains to show that the conditions (iii) and (iv) hold at
points in $\what X\sm(\what V_1\cup\cdots\cup\what
V_k)\subset\what W_0$. But $\alpha\equiv 0$ on $\what W_0$, so the
condition~(iii) holds for $R>R_0$ (or even $R>0$) and we have
$\eta_-\left(\dist_{\hat g}(\cdot,Z)\right)\leq\alpha$. Moreover,
on $\what W_0$, we have
\[
\dist_{\hat g}(\cdot,Z)\geq\dist_{\hat g}(\cdot,\what
Y)\geq\dist_{\hat g}(\what V_0,\what X\sm\what
U_0)=\dist_g(V_0,X\sm U_0).
\]
and hence
\[
\alpha=0=\eta_+(\dist_{\hat g}(\cdot,Z)).
\]
Thus $\alpha$ satisfies the conditions (i)--(iv) (for $R>R_0$) on
the entire covering space $\what X$.
\end{pf*}

\section{An $\strplshclass^\infty(g,q)$ \fn near the singular set}\label{q-plsh
fn near the singular set sect}

For the proof of Theorem~\ref{main theorem from intro},
Proposition~\ref{q-cvx on nbd off sing set prop} allows one to
construct a \fn which, in a \nbd of $\what Y$, is \strg $q$-\cvx
away from $\sing{\what Y}$. In order to modify the \fn to be \strg
$q$-\cvx near the singular set, we will apply the following fact;
the proof of which is the first goal of this section:
\begin{prop}\label{q-plsh fn near sing set prop}
Let $(X,g)$ be a \con reduced Hermitian \cpx space, let $q$ be a
positive integer, let $Y$ be a \cpt \anal subset of~$X$, let $S$
be a \cpt \anal subset of dimension~$<q$ in~$X$, let $\Omega$ be a
\rel \cpt \nbd of~$S$ in~$X$, let $\seq D\nu_{\nu=1}^m$ be \rel
\cpt open subsets of~$X$, and let $\rho_\nu$ be a \cinf \str \plsh
\fn on a \nbd of $\overline D_\nu$ for each $\nu=1,\dots,m$. Then
there is a \rel \cpt \nbdns~$\Theta$ of $S$ in $\Omega$, a
constant $R>0$, and, for every $\epsilon>0$, an associated
$\delta>0$ \stns, for every \con covering space
$\Upsilon\colon\what X\to X$ and for every \anal set $Z$ equal to
a union of \ircomps of $\what Y\equiv\Upsilon\inv(Y)$, there
exists a \cinf \fn $\alpha$ on $\what X$ with the following
properties relative to the sets $\what\Omega=\Upsilon\inv(\Omega)$
and $\what\Theta=\Upsilon\inv(\Theta)$ and the Hermitian metric
$\hat g=\Upsilon^*g$:
\begin{enumerate}

\item[(i)] On $\what X$, $0\leq\alpha\leq R$;

\item[(ii)] We have $\supp\alpha\subset\what\Omega\sm\deltanbd
Z\delta$ and $\alpha>0$ on $\what\Theta\sm\deltanbd Z\epsilon$;

\item[(iii)] On $\what\Theta\cup\deltanbd Z\delta$, $\alpha$ is of
class $\plshclass^\infty(\hat g,q)$;

\item[(iv)] On the set $\setof{p\in\what\Theta}{\alpha(p)>0}$,
$\alpha$ is of class $\strplshclass^\infty(\hat g,q)$; and

\item[(v)] For each $\nu=1,\dots,m$,
$\alpha+R\rho_\nu\circ\Upsilon$ is \cinf \str \plsh on
$\Upsilon\inv(D_\nu)$ and $R\alpha-\rho_\nu\circ\Upsilon$ is of
class $\strplshclass^\infty(\hat g,q)$ on a \nbd of
$\Upsilon\inv(D_\nu)\cap\what\Theta\sm\deltanbd Z\epsilon$.

\end{enumerate}
\end{prop}
\begin{rmk}
Note that $\Theta$ and $R$ do not depend on the choice of
$\epsilon$ and $\delta$.
\end{rmk}
\begin{pf}
We may assume without loss of generality that $\Omega\Subset
D_1\cup\cdots\cup D_m$.  According to Theorem~\ref{Qplsh on nbd of
lower dimension or noncpt q dim thm}, there exists a nonnegative
\cinf \fn $\beta$ on $X$ \st $\supp\beta\subset\Omega$ and \stns,
on a \nbd of the closure of some \rel \cpt \nbdns~$\Theta$ of $S$
in $\Omega$, we have $\beta>2$ and $\beta$ is of class
$\strplshclass^\infty(g,q)$. We may also choose a \cinf \fn
$\chi\colon\R\to[0,\infty)$ \st $\chi\equiv 0$ on $(-\infty,0]$,
$\chi'>0$ and $\chi''\geq 0$ on $(0,\infty)$, and $\chi'(1)\geq
1$. Choosing a constant $R_0\gg 1$, we get $\beta\leq R_0$,
$\chi(\beta)\leq R_0$, and $\chi'(\beta)\leq R_0$ on $X$. Choosing
$R_1\gg 0$, we get, for each $\nu=1,\dots,m$, $R_1\beta-\rho_\nu$
is of class $\strplshclass^\infty(g,q)$ on $D_\nu\cap\Theta$ and
$R_0\beta+R_1\rho_\nu$ is \cinf \str \plsh on $D_\nu$.

Applying Proposition~\ref{Quasi-plsh fn near the singular set
prop}, we may choose (independently of the choices made in the
previous paragraph) a \rel \cpt \nbdns~$\Lambda$ of $Y$ in $X$, a
pair of nondecreasing \cont maps
\[
\eta_-\colon[0,\infty]\to[-\infty,0]\qquad\text{and}\qquad\eta_+\colon
[0,\infty]\to [-\infty,0]
\]
with \(\eta_-\leq\eta_+\), \(\eta_+(0)=-\infty\), and
\(\eta_->-\infty\) on $(0,\infty]$, and a constant $R_2>0$ \stns,
for every \con covering space $\Upsilon\colon\what X\to X$ and for
every \anal set $Z$ equal to a union of \ircomps of $\what
Y\equiv\Upsilon\inv(Y)$, there exists a \cont \fn
$\gamma\colon\what X\to [-\infty,0]$ with the following
properties:
\begin{enumerate}
\item[(\ref{q-plsh fn near sing set prop}.1)] The \fn $\gamma$ is
of class $\cal Q^\infty_Z$ on $\what X$;

\item[(\ref{q-plsh fn near sing set prop}.2)] We have
$\supp\gamma\subset\hat\Lambda\equiv\Upsilon\inv(\Lambda)$;

\item[(\ref{q-plsh fn near sing set prop}.3)] For each
$\nu=1,\dots,m$, $\gamma+R_2\cdot\rho_\nu\circ\Upsilon$ is of
class $\strplshclass^\infty_Z$ on $\Upsilon\inv(D_\nu)$;

\item[(\ref{q-plsh fn near sing set prop}.4)] For $\hat
g=\Upsilon^*g$, we have
\[
\eta_-(\dist_{\hat g}(x,Z))\leq\gamma(x)\leq\eta_+(\dist_{\hat
g}(x,Z))\qquad\forall\, x\in\what X.
\]
\end{enumerate}

We may now fix a constant $R>\max(R_0,2R_1)$ and, given
$\epsilon>0$, we may choose a constant $\mu>0$ so small that $2\mu
R_1R_2<2\mu R_0R_1R_2<1$, $R_1+\mu R_0R_2<R$, and
$\mu\cdot\eta_-(\epsilon)>-1$; and we may then choose a constant
$\delta\in(0,\epsilon)$ so small that
$\mu\cdot\eta_+(\delta)<-R_0$. We will show that, for
$\Upsilon\colon\what X\to X$ a \con covering space, $Z$ a union of
\ircomps of $\what Y=\Upsilon\inv(Y)$, and $\gamma\colon\what
X\to[-\infty,0]$ satisfying (\ref{q-plsh fn near sing set
prop}.1)--(\ref{q-plsh fn near sing set prop}.4), the \fn $\alpha$
given by $\alpha\equiv 0$ on $Z$ and
$\alpha=\chi(\beta\circ\Upsilon+\mu\gamma)$ on $\what X\sm Z$ has
the required properties (i)--(v). For this, we set
$\what\Omega=\Upsilon\inv(\Omega)$,
$\what\Theta=\Upsilon\inv(\Theta)$,
$\hat\beta=\beta\circ\Upsilon$, and, for each $\nu=1,\dots,m$,
$\what D_\nu=\Upsilon\inv(D_\nu)$ and
$\hat\rho_\nu=\rho_\nu\circ\Upsilon$ on $\what D_\nu$.

The choice of $\chi$ and $R_0$ guarantee that $\alpha$ is of class
\cinf and $0\leq\alpha\leq\chi(\hat\beta)\leq R$ on $\what X$; as
in~(i). We also have $\hat\beta+\mu\gamma\leq 0$ on
$\hat\beta\inv(0)\cup\deltanbd Z\delta$ and
$\hat\beta+\mu\gamma\geq\hat\beta+\mu\eta_-(\epsilon)>2-1>0$ on
$\what\Theta\sm\deltanbd Z\epsilon$, so the condition~(ii) holds.
For each $\nu=1,\dots,m$, we have, on the set
$\what\Theta\cap\what D_\nu\sm Z$,
\[
\hat\beta+\mu\gamma=R_1\inv\cdot(R_1\hat\beta-\hat\rho_\nu)+\mu\cdot(\gamma+R_2\hat\rho_\nu)
+(R_1\inv-\mu R_2)\cdot\hat\rho_\nu.
\]
Since $\chi',\chi''\geq 0$ and $\alpha\equiv 0$ on $\deltanbd
Z\delta$, Lemma~\ref{convex fn comp lemma} implies that $\alpha$
satisfies the condition~(iii). Furthermore, since $\chi'(t)>0$
precisely when $\chi(t)>0$ (i.e. when $t>0$), the condition~(iv)
holds.

It remains to verify the condition~(v). Given an index $\nu$ with
$1\leq\nu\leq m$ and a point $p\in\what D_\nu$, we may choose a
proper local \holo model $(U,\Phi,U')$ \st $p\in U\subset\what
D_\nu$ and \stns, for some choice of \cinf \str \plsh \fns $\tau$
and $\rho$ on $U'$ and some \fn $\theta$ of class
$\strplshclass^\infty_{\Phi(Z\cap U)}$ on $U'$, we have
$R_0\hat\beta+R_1\hat\rho_\nu=\tau\circ\Phi$,
$\hat\rho_\nu=\rho\circ\Phi$, and
$\gamma+R_2\cdot\hat\rho_\nu=\theta\circ\Phi$ on $U$. In
particular, we get the extensions
\[
\beta_0\equiv R_0\inv\cdot(\tau-R_1\rho)\quad\text{and}
\quad\gamma_0\equiv\theta-R_2\cdot\rho
\]
of $\hat\beta\restrict U$ and $\gamma\restrict U$, respectively.
If $p\in\what\Theta\cap\what D_\nu$, then we may also choose the
local model so that $U\subset\what\Theta$ and so that, for some
\fn $\omega\in\cinf(U')$, we have
$R_1\hat\beta-\hat\rho_\nu=\omega\circ\Phi$ on $U$ and the $\hat
g$-trace of the restriction of $\Phi^*\lev\omega$ to any
$q$-dimensional subspace of $\ztanq 1_xU$ is positive for each
$x\in U$. We then get a second extension \(\beta_1\equiv
R_1\inv\cdot(\omega+\rho) \) of $\hat\beta\restrict U$.

We also have a \cinf extension $\vphi$ of the \fn
$(\alpha+R\hat\rho_\nu)\restrict U$ given by
\[
\vphi=\left\{
\begin{aligned}
R\rho&\qquad\text{on }\Phi(Z\cap U)\\
\chi(\beta_0+\mu\gamma_0)+R\rho&\qquad\text{on }U'\sm\Phi(Z\cap U)
\end{aligned}
\right.
\]
In particular, $\vphi=R\rho$ near $\Phi(Z\cap U)$ and hence
$\vphi$ is \str \plsh near $\Phi(Z\cap U)$. For each point
$x\in\Phi(U\sm Z)$ and each nonzero tangent vector $v\in
T^{1,0}_xU'$, we have
\[
\lev\vphi(v,v)\geq\chi'\left[\beta_0(x)+\mu\gamma_0(x)\right]\cdot\lev{\beta_0+\mu\gamma_0}(v,v)+R\cdot\lev\rho(v,v).
\]
Hence $\lev\vphi(v,v)>0$ if $\lev{\beta_0+\mu\gamma_0}(v,v)\geq
0$. If $\lev{\beta_0+\mu\gamma_0}(v,v)<0$, then, since
$0\leq\chi'(\beta_0(x)+\mu\gamma_0(x))\leq\chi'(\beta_0(x))\leq
R_0$, we get
\begin{align*}
\lev\vphi(v,v)&\geq
R_0\cdot\lev{\beta_0+\mu\gamma_0}(v,v)+R\cdot\lev\rho(v,v)\\
&=R_0\cdot\lev{R_0\inv\cdot(\tau-R_1\rho)+
\mu\cdot(\theta-R_2\cdot\rho)}(v,v)+R\cdot\lev\rho(v,v)\\
&=\lev\tau(v,v)+\mu R_0\cdot\lev\theta(v,v)+(R-R_1-\mu
R_0R_2)\cdot\lev\rho(v,v)\\
&>0.
\end{align*}
It follows that $\vphi$ is \str \plsh on a \nbd of $\Phi(U)$.

If $p\in\what\Theta\cap D_\nu$, $\psi$ is the \cinf extension of
the \fn $(R\alpha-\hat\rho_\nu)\restrict{U\sm Z}$ to
$U'\sm\Phi(Z\cap U)$ given by $\psi\equiv
R\cdot\chi(\beta_1+\mu\gamma_0)-\rho$, $y\in U\sm Z$, $x=\Phi(y)$,
$e_1,\dots,e_q\in\ztanq 1_yX$ are orthonormal, and $v_j=\Phi_*e_j$
for $j=1,\dots,q$, then
\begin{align*}
\sum_{j=1}^q\lev\psi(v_j,v_j)&\geq
R\cdot\chi'\left(\hat\beta(y)+\mu\gamma(y)\right)\cdot\sum_{j=1}^q\lev{\beta_1+\mu\gamma_0}(v_j,v_j)-\sum_{j=1}^q\lev\rho(v_j,v_j)\\
&=RR_1\inv\cdot\chi'\left(\hat\beta(y)+\mu\gamma(y)\right)\cdot\sum_{j=1}^q\lev\omega(v_j,v_j)\\
&\qquad+\mu
R\cdot\chi'\left(\hat\beta(y)+\mu\gamma(y)\right)\cdot\sum_{j=1}^q\lev\theta(v_j,v_j)\\
&\qquad+\left[R(R_1\inv-\mu
R_2)\cdot\chi'\left(\hat\beta(y)+\mu\gamma(y)\right)-1\right]\cdot\sum_{j=1}^q\lev\rho(v_j,v_j)
\end{align*}
We have $R_1\inv-\mu R_2=R_1\inv(1-\mu R_1R_2)>(2R_1)\inv$. If
$y\in U\sm\deltanbd Z\epsilon$, then
$\hat\beta(y)+\mu\gamma(y)>2+\mu\cdot\eta_-(\epsilon)>1$ and hence
\[
R(R_1\inv-\mu
R_2)\cdot\chi'(\hat\beta(y)+\mu\gamma(y))-1>(R/(2R_1))-1>0.
\]
Thus $\sum\lev\psi(v_j,v_j)>0$ for any choice of $y$ in a small
\nbd of $U\sm\deltanbd Z\epsilon$ and hence $\alpha$ and $R$
satisfy the condition~(v).
\end{pf}

Combining Proposition~\ref{q-cvx on nbd off sing set prop} and
Proposition~\ref{q-plsh fn near sing set prop}, we get the
following:
\begin{prop}\label{q-cvx away from cpt not exhaustive prop}
Let $(X,g)$ be a \con reduced Hermitian \cpx space; let $q$ be a
positive integer; let $Y$ be a \cpt \anal set of dimension~$\leq
q$; let~$C$~and~$S$ be \anal subsets of~$Y$ \st $C$ is a union of
\ircomps of $Y$, $\dim S<q$, $S$ contains $\sing Y$ as well as
every \ircomp of $Y$ of dimension~$<q$, and $Y\sm (C\cup S)$ is
Stein; let $\seq D\nu_{\nu=1}^m$ be \rel \cpt open subsets of~$X$;
let $\rho_\nu$ be a \cinf \str \plsh \fn on a \nbd of $\overline
D_\nu$ for each $\nu=1,\dots,m$; and let $\epsilon>0$. Then, for
every choice of constants $\delta_1$,~$\delta_2$,~$\delta_3$,~and
$\delta_4$ with
$\epsilon\gg\delta_4\gg\delta_3\gg\delta_2\gg\delta_1>0$ (i.e.~one
must choose $\delta_4$ sufficiently small relative to~$\epsilon$,
$\delta_3$ sufficiently small relative to $\delta_4$, and so on),
there exists a \nbd $\Omega$ of $Y$ in $X$ \stns, for every \con
covering space $\Upsilon\colon\what X\to X$, for every \anal set
$Z$ which is equal to a union of \ircomps of $\what
Y\equiv\Upsilon\inv(Y)$ and which contains $\what
C=\Upsilon\inv(C)$ as well as every $q$-dimensional \cpt \ircomp
of $\what Y$ (i.e. $Z$ contains every \cpt \ircomp of $\what Y$
not contained in $\what S=\Upsilon\inv(S)$), and for every
positive \cont \fn $\theta$ on $\what X$, there exists a
nonnegative \cinf \fn $\alpha$ on $\what X$ with the following
properties relative to the Hermitian metric $\hat g=\Upsilon^*g$,
the sets $\what\Omega=\Upsilon\inv(\Omega)$ and $\what
D_\nu=\Upsilon\inv(D_\nu)$ for $\nu=1,\dots,m$, and the \fns
$\hat\rho_\nu=\rho_\nu\circ\Upsilon\colon\what D_\nu\to\R$ for
$\nu=1,\dots,m$:
\begin{enumerate}

\item[(i)] On $\deltanbd Z{\delta_1}\cap\what\Omega$,
$\alpha\equiv 0$;

\item[(ii)] On $\left(\what\Omega\sm\deltanbd
Z{\delta_2}\right)\cup\left(\what\Omega\sm\left[\deltanbd{\what
S}{\delta_2}\cup\deltanbd Z{\delta_1}\right]\right)$, $\alpha>0$;

\item[(iii)] On $\deltanbd{\what S}{\delta_3}\cap\what\Omega$,
$\alpha$ is of class $\plshclass^\infty(\hat g,q)$;

\item[(iv)] On $\setof{p\in\deltanbd {\what
S}{\delta_3}\cap\what\Omega}{\alpha(p)>0}$, $\alpha$ is of class
$\strplshclass^\infty(\hat g,q)$;

\item[(v)] For each $\nu=1,\dots,m$, the \fn $\alpha-\hat\rho_\nu$
is of class $\strplshclass^\infty(\hat g,q)$ on a \nbd of
$\deltanbd{\what S}{\delta_3}\cap\what
D_\nu\cap\what\Omega\sm\deltanbd Z{\delta_2}$;

\item[(vi)] For each $\nu=1,\dots,m$, the \fn
$\alpha-\hat\rho_\nu$ is \cinf \strg $q$-\cvx on a \nbd of $\what
D_\nu\cap\what\Omega\sm\left[\deltanbd{\what
S}{\delta_2}\cup\deltanbd Z{\delta_1}\right]$;

\item[(vii)] On $\what\Omega\sm\left[\deltanbd{\what
S}{\delta_4}\cup\deltanbd Z{\delta_1}\right]$, $\alpha>\theta$;
and

\item[(viii)] On $\setof{p\in\what\Omega}{\alpha(p)>0}$, $\alpha$
is \cinf \strg $q$-\cvxns.

\end{enumerate}
\end{prop}
\begin{pf}
We may assume without loss of generality that $Y\subset
\Theta_0\Subset D_1\cup\cdots\cup D_m$ for some open set
$\Theta_0$. We will also assume that $Y\sm (C\cup
S)\neq\emptyset$, since the proof for $Y\subset C\cup S$ is
similar (but easier). Applying Proposition~\ref{q-plsh fn near
sing set prop}, we get a \rel \cpt \nbdns~$\Theta_1$ of $S$ in
$\Theta_0$, a constant $R_0>1$, and, for every $\eta>0$, an
associated $\delta>0$ \stns, for every \con covering space
$\Upsilon\colon\what X\to X$ and for every \anal set $Z$ equal to
a union of \ircomps of $\what Y\equiv\Upsilon\inv(Y)$, there
exists a nonnegative \cinf \fn $\beta$ on $\what X$ with the
following properties relative to the Hermitian metric $\hat
g=\Upsilon^*g$, the sets $\what\Theta_j=\Upsilon\inv(\Theta_j)$
for $j=0,1$ and $\what D_\nu=\Upsilon\inv(D_\nu)$ for
$\nu=1,\dots,m$, and the \fns
$\hat\rho_\nu=\rho_\nu\circ\Upsilon\colon\what D_\nu\to\R$ for
$\nu=1,\dots,m$:
\begin{enumerate}

\item[(\ref{q-cvx away from cpt not exhaustive prop}.1)] We have
$\supp\beta\subset\what\Theta_0\sm\deltanbd Z\delta$ and $\beta>0$
on $\what\Theta_1\sm\deltanbd Z\eta$;

\item[(\ref{q-cvx away from cpt not exhaustive prop}.2)] On
$\what\Theta_1\cup\deltanbd Z\delta$, $\beta$ is of class
$\plshclass^\infty(\hat g,q)$;

\item[(\ref{q-cvx away from cpt not exhaustive prop}.3)] On the
set $\setof{p\in\what\Theta_1}{\beta(p)>0}$, $\beta$ is of class
$\strplshclass^\infty(\hat g,q)$; and

\item[(\ref{q-cvx away from cpt not exhaustive prop}.4)] For each
$\nu=1,\dots,m$, $\beta+R_0\hat\rho_\nu$ is \cinf \str \plsh on
$\what D_\nu$ and $\beta-\hat\rho_\nu$ is of class
$\strplshclass^\infty(\hat g,q)$ on a \nbd of $\what
D_\nu\cap\what\Theta_1\sm\deltanbd Z\eta$ (in the notation of
Proposition~\ref{q-plsh fn near sing set prop}, $\beta=R\alpha$
and $R_0=R^2$).

\end{enumerate}

We may now choose constants $\xi$, $\delta_2$, $\delta_3$,
$\delta_4$ and open sets $\Theta_2$, $\Omega_1,\dots,\Omega_4$
with the following properties:
\begin{enumerate}
\item[(\ref{q-cvx away from cpt not exhaustive prop}.5)] We have
$0<2\xi<\delta_2<\delta_3<\delta_4<\min\left[\epsilon,\dist(Y,X\sm\Theta_0),\dist(S,X\sm\Theta_1)\right]$;

\item[(\ref{q-cvx away from cpt not exhaustive prop}.6)] We have
$\Omega_4\Subset\cdots\Subset\Omega_1\Subset\Theta_0\sm(C\cup S)$;

\item[(\ref{q-cvx away from cpt not exhaustive prop}.7)] For each
\concomp $A$ of $Y\sm(C\cup S)$ and each $j=1,\dots,4$, the set
$\Omega_j\cap A$ is nonempty and \conns, $\overline{\Omega_j\cap
A}=\overline\Omega_j\cap A$, and
\[
\pi_1(\Omega_j\cap
A)\twoheadrightarrow\text{im}\,\left[\pi_1(A)\to\pi_1(\bar
A)\right];
\]

\item[(\ref{q-cvx away from cpt not exhaustive prop}.8)]
\(Y\sm(C\cup S)\subset\Theta_2\subset\Theta_0\sm(C\cup S)\),
\(\Theta_2\sm\deltanbd S{\delta_4}\subset\Omega_4\),
\(\Omega_3\cap\deltanbd{C\cup S}{2\delta_3}=\emptyset\),
\(\Theta_2\sm\deltanbd S{\delta_2}\subset\Omega_2\), and
\(\Omega_1\cap\deltanbd{C\cup S}{2\xi}=\emptyset\);

\item[(\ref{q-cvx away from cpt not exhaustive prop}.9)] For any
\con covering space $\what X$ and any \anal set $Z\subset\what X$
as above, we may form a \fn $\beta$ with the properties
(\ref{q-cvx away from cpt not exhaustive prop}.1)--(\ref{q-cvx
away from cpt not exhaustive prop}.4) for $\eta=\delta_2$ and
$\delta=\xi$.

\end{enumerate}
Note that we obtain the above sets and constants by choosing them
in the order: $\delta_4$, $\Omega_4$, $\Omega_3$, $\delta_3$,
$\delta_2$, $\Omega_2$, $\Omega_1$, $\xi$, $\Theta_2$. To get the
surjection of fundamental groups in (\ref{q-cvx away from cpt not
exhaustive prop}.7), we choose $\Omega_4$ sufficiently large and
apply standard facts (see, for example, Lemma~2
of~\cite{Fraboni-Covering cpx mfld}). To get $\xi$ to
satisfy~(\ref{q-cvx away from cpt not exhaustive prop}.9), we need
only choose $\eta$ so that $0<\eta\leq\delta_2$ and
\(\Omega_1\cap\deltanbd{C\cup S}{\eta}=\emptyset\), and then
choose $\xi$ to be an associated $\delta$ as in (\ref{q-cvx away
from cpt not exhaustive prop}.1)--(\ref{q-cvx away from cpt not
exhaustive prop}.4) with $0<2\xi<\eta$. To get the condition
(\ref{q-cvx away from cpt not exhaustive prop}.8), we first choose
the constants $\xi$, $\delta_2$, $\delta_3$, $\delta_4$ and the
sets $\Omega_i$ for $i=1,\dots,4$ so that (\ref{q-cvx away from
cpt not exhaustive prop}.8) holds with $Y\sm(C\cup S)$ in place of
$\Theta_2$, and we then choose the \nbd $\Theta_2$ sufficiently
small.

For each \concomp $Y'$ of $Y\sm(C\cup S)$, we may form a \con \nbd
$X'$ in $\Theta_2$ with $\overline{X'}\subset\Theta_2\cup S$ and
$\pi_1(Y')\twoheadrightarrow\text{im}\,\left[\pi_1(X')\to\pi_1(X)\right]$.
We may also form a \nbd $\Xi'$ of $Y'$ in $X'$ and open sets
$\Omega_1',\dots,\Omega_4'$ with
\[
X'\Supset\Omega_1'\Supset\cdots\Supset\Omega_4'\qquad\text{and}\qquad
\Omega_j'\cap\Xi'=\Omega_j\cap\Xi'\text{ for }j=1,\dots,4.
\]
By applying Proposition~\ref{q-cvx on nbd off sing set prop}, we
get a corresponding \nbd as in part~(b) (of Proposition~\ref{q-cvx
on nbd off sing set prop}) which we may take to be contained in
$\Xi'$ (note that $X'$ is chosen so that each \comp $\what X'$ of
the lifting in any covering space of $X$ will contain exactly one
\comp $\what Y'$ of the lifting of $Y'$, so the proposition may be
applied to the covering space $\what X'\to\what Y'$ whenever this
restricted covering is infinite).

Taking the union of the (finitely many) resulting \nbds of all of
the \concomps of $Y\sm(C\cup S)$ (which we may take to be
disjoint), we get a \nbd~$\Theta_3$ of $Y\sm(C\cup S)$ in
$\Theta_2$ and an open subset $B$ of $\Theta_3$ \st
$\overline\Theta_3\subset\Theta_2\cup S$, each
\concompns~$\Theta'$ of $\Theta_3$ meets (hence contains) exactly
one \concompns~$Y'$ of $Y\sm(C\cup S)$ and satisfies
\[
\pi_1(Y')\twoheadrightarrow\text{im}\,\left[\pi_1(\Theta')\to\pi_1(X)\right],
\]
and, for every connected covering space $\Upsilon\colon\what X\to
X$ and every positive \cont \fn $\theta$ on $\what X$, there is a
nonnegative \cinf \fn $\gamma$ on the union $\Gamma_0$ of all of
the \concomps of $\what\Theta_3=\Upsilon\inv(\Theta_3)$ which meet
(hence contain) a \concomp of $\Upsilon\inv\left(Y\sm(C\cup
S)\right)$ which is not \rel \cpt in $\what X$ \stns, if
$\what\Omega_j=\Upsilon\inv(\Omega_j)$ for $j=1,2,3,4$, then
\begin{enumerate}
    \item[(\ref{q-cvx away from cpt not exhaustive prop}.10)] On $\Gamma_0\sm\what\Omega_1$, $\gamma\equiv 0$;
    \item[(\ref{q-cvx away from cpt not exhaustive prop}.11)] On $\what\Omega_2\cap\Gamma_0$, $\gamma>0$;
    \item[(\ref{q-cvx away from cpt not exhaustive prop}.12)] There is a nonnegative \fn $\gamma_0\in\plshclass^\infty(g,q)(\Theta_3\sm\overline\Omega_4)$
    \st $\gamma=\gamma_0\circ\Upsilon$ on
    $\Gamma_0\sm\what\Omega_3$ and $\gamma_0$ is of class $\strplshclass^\infty(g,q)$
    on $\setof{x\in\Theta_3\sm\overline\Omega_4}{\gamma_0(x)>0}
    \supset\Omega_2\cap\Theta_3\sm\overline\Omega_4$;
    \item[(\ref{q-cvx away from cpt not exhaustive prop}.13)]
    We have
    \(B=(\Omega_3\cap\Theta_3)\cup\setof{x\in\Theta_3\sm\Omega_3}{\gamma_0(x)>0}\supset\Omega_2\cap\Theta_3\)
    and, for $\what B=\Upsilon\inv(B)$, we have
    \[
    \what B\cap\Gamma_0=\setof{x\in\Gamma_0}{\gamma(x)>0}
    \]
    and, for any \cinf \fn $\tau$ with \cpt support in~$B$, the \fn
$R\cdot\gamma+\tau\circ\Upsilon$ will be
    \cinf \strg $q$-\cvx on $\what B\cap\Gamma_0$
    for every sufficiently large positive constant $R$; and
    \item[(\ref{q-cvx away from cpt not exhaustive prop}.14)] On $\what\Omega_4\cap\Gamma_0$,
    $\gamma>\theta$.
\end{enumerate}

Finally, we may choose a constant $\delta_1$ \st $0<\delta_1<\xi$,
$\dist(\Theta_3\sm\deltanbd S{\xi},C)>\delta_1$, and
$\dist\left(Y\sm(\deltanbd S{\xi}\cup
C),X\sm\Theta_3\right)>\delta_1$.

We will now show that any \rel \cpt \nbd $\Omega$ of $Y$ in the
\nbd $\Omega_0\equiv\deltanbd
Y{\delta_1}\cap\left[\Theta_3\cup\deltanbd {C\cup
S}{\delta_1}\right]$ and the constants $\delta_1$, $\delta_2$,
$\delta_3$, and $\delta_4$ have the required properties. For this,
suppose $\Upsilon\colon\what X\to X$ is a \con covering space,
$\theta$ is a positive \cont \fn on $\what X$, $Z$ is a union of
\ircomps of $\what Y=\Upsilon\inv(Y)$ which contains $\what
C=\Upsilon\inv(C)$ as well as every $q$-dimensional \cpt \ircomp
of $\what Y$, $\hat g=\Upsilon^*g$,
$\what\Theta_i=\Upsilon\inv(\Theta_i)$ for $i=0,1,2,3$,
$\what\Omega_i=\Upsilon\inv(\Omega_i)$ for $i=0,1,2,3,4$,
$\what\Omega=\Upsilon\inv(\Omega)$, and
$\hat\rho_\nu=\rho_\nu\circ\Upsilon$ on $\what
D_\nu\equiv\Upsilon\inv(D_\nu)$ for $\nu=1,\dots,m$. We will
assume that the covering is infinite, since the proof for a finite
covering is similar (but easier). We may also assume that $\theta$
is an \exh \fn (since we may replace $\theta$ by an \exh \fn which
is greater than $\theta$). In the above notation, we may form a
nonnegative \cinf \fn $\beta$ on $\what X$ satisfying the
conditions (\ref{q-cvx away from cpt not exhaustive
prop}.1)--(\ref{q-cvx away from cpt not exhaustive prop}.4) with
$\eta=\delta_2$ and $\delta=\xi$, and we may form a nonnegative
\cinf \fn $\gamma$ on $\Gamma_0\subset\what\Theta_3$ satisfying
the conditions (\ref{q-cvx away from cpt not exhaustive
prop}.10)--(\ref{q-cvx away from cpt not exhaustive prop}.14).
Observe that, for $B$ and $\what B=\Upsilon\inv(B)$ as in
(\ref{q-cvx away from cpt not exhaustive prop}.13), we have
\(\Omega\cap\Theta_3\sm\deltanbd S{\delta_2}\Subset
U\Subset\Omega_2\cap\Theta_3\subset B \) for some open set~$U$.
For, if $p\in\Omega\cap\Theta_3\sm\deltanbd S{\delta_2}$, then
$\dist(p,x)<\delta_1$ for some point $x\in Y$ and we get
$\dist(x,S)>\delta_2-\delta_1>\xi$. Since
$\delta_1<\dist(\Theta_3\sm\deltanbd S\xi,C)$, we also have
$x\not\in C$. Thus
\[
\Omega\cap\Theta_3\sm\deltanbd
S{\delta_2}\subset\deltanbd{Y\sm\left[\deltanbd S\xi\cup
C\right]}{\delta_1}\sm\deltanbd
S{\delta_2}\Subset\Theta_3\sm\deltanbd
S{\delta_2}\subset\Omega_2\cap\Theta_3.
\]
Setting $\what U=\Upsilon\inv(U)$, we see that, by replacing
$\gamma$ with a large multiple, we may assume that
$\gamma-(R_0+1)\cdot\hat\rho_\nu$ is \cinf \strg $q$-\cvx on
$\what D_\nu\cap\what U\cap\Gamma_0$ for $\nu=1,\dots,m$.

Let $\Gamma\subset\Gamma_0$ be the union of all of the \concomps
of $\what\Theta_3$ which meet (hence contain) a \concomp of $\what
Y\sm(Z\cup\what S)$. Observe that we get a well-defined
nonnegative \cinf \fn $\zeta$ on $\what\Omega_0$ by setting
$\zeta=\beta+\gamma$ on $\Gamma\cap\what\Omega_0$ and $\beta$
elsewhere in $\what\Omega_0$. For
\[
\what\Omega_0\cap\partial\Gamma\subset\deltanbd{\what
Y}{\delta_1}\cap\left[\what\Theta_3\cup\deltanbd{\what C\cup\what
S}{\delta_1}\right]\cap\partial\what\Theta_3=\deltanbd{\what
C\cup\what S}{\delta_1}\cap\partial\what\Theta_3.
\]
Since $\gamma\equiv 0$ on
$\Gamma_0\sm\what\Omega_1\supset\deltanbd{\what C\cup\what
S}{\delta_1}\cap\Gamma_0$, we see that $\zeta$ is $\cinf$. Cutting
off, we get a nonnegative \cinf \fn $\alpha$ on $\what X$ with
$\alpha=\zeta$ on $\what\Omega$.

In order to verify the properties (i)--(viii), we first observe
that we have the following:
\[
\what\Omega\subset\deltanbd{\what
Y}{\delta_1}\subset\deltanbd{\what S}{2\xi}\cup\Gamma\cup\deltanbd
Z{\delta_1}\qquad\text{and}\qquad\Gamma\cap\deltanbd
Z{\delta_1}\subset\deltanbd{\what S}{2\xi}.
\]
For, given $x\in\deltanbd{\what Y}{\delta_1}\sm\deltanbd{\what
S}{2\xi}$, we may fix a point $y\in\what Y$ with
$\dist(x,y)<\delta_1$. Hence \( \dist\left(y,\what
S\right)>2\xi-\delta_1>\xi\). If $y\in\what C\subset Z$, then,
since $\dist\left[\what\Theta_3\sm\deltanbd{\what S}\xi,\what
C\right]>\delta_1$, we get $x\in\deltanbd
Z{\delta_1}\sm\what\Theta_3\subset\deltanbd Z{\delta_1}\sm\Gamma$.
If $y\notin\what C$, then, since $\dist\left[\what
Y\sm\left(\deltanbd {\what S}{\xi}\cup \what C\right),\what
X\sm\what\Theta_3\right]>\delta_1$, a piecewise \cinf path from
$y$ to $x$ of length~$<\delta_1$ must lie entirely in the \concomp
$V$ of $\what\Theta_3$ containing~$y$. Thus either $y\in Z$, in
which case $x\in\deltanbd Z{\delta_1}$ and  $x\in
V\subset\what\Theta_3\sm\Gamma$; or $y\notin Z$, in which case
$x\in V\subset\Gamma$.

\textbf{Verification of (i).} The condition~(i) holds because
$\gamma\equiv 0$ on
$\Gamma\sm\what\Omega_1\supset\Gamma\cap\deltanbd{\what C\cup\what
S}{2\xi}\supset\Gamma\cap\deltanbd Z{\delta_1}$ and $\beta\equiv
0$ on $\deltanbd Z\xi\supset\deltanbd Z{\delta_1}$.

\textbf{Verification of (ii).} We have $\alpha\geq\beta>0$ on
$\what\Theta_1\cap\what\Omega\sm\deltanbd Z{\delta_2}$ and
$\alpha\geq\gamma>0$ on
$\what\Omega_2\cap\Gamma\cap\what\Omega\supset\Gamma\cap\what\Omega\sm\deltanbd{\what
S}{\delta_2}\supset\Gamma\cap\what\Omega\sm\deltanbd{\what
S}{\delta_4}\supset\Gamma\cap\what\Omega\sm\what\Theta_1$. On the
other hand, we have
\[
\what\Omega\sm\left(\Gamma\cup\what\Theta_1\right)\subset
\what\Omega\sm\left(\Gamma\cup\deltanbd{\what
S}{\delta_2}\right)\subset
\what\Omega\sm\left(\Gamma\cup\deltanbd{\what
S}{2\xi}\right)\subset\deltanbd Z{\delta_1}.
\]
That~(ii) holds now follows.

\textbf{Verification of (iii).} From (\ref{q-cvx away from cpt not
exhaustive prop}.2), we see that $\beta$ is of class
$\plshclass^\infty(\hat g,q)$ on
$\what\Theta_1\supset\deltanbd{\what S}{\delta_3}$; and from
(\ref{q-cvx away from cpt not exhaustive prop}.8) and (\ref{q-cvx
away from cpt not exhaustive prop}.12), we see that $\gamma$ is of
class $\plshclass^\infty(\hat g,q)$ on
$\Gamma\sm\overline{\what\Omega_3}\supset\Gamma\cap\deltanbd{\what
S}{2\delta_3}$. It follows that the condition~(iii) holds.

\textbf{Verification of (iv).} Suppose $p\in\deltanbd{\what
S}{\delta_3}\cap\what\Omega$ with $\alpha(p)=\zeta(p)>0$. If
$\beta(p)>0$, then $\beta$ is of class $\strplshclass^\infty(\hat
g,q)$ near $p$ (by (\ref{q-cvx away from cpt not exhaustive
prop}.3)) and hence, since $\gamma$ is of class
$\plshclass^\infty(\hat g,q)$ on $\Gamma\cap\deltanbd{\what
S}{\delta_3}$, we see that $\alpha$ is of class
$\strplshclass^\infty(\hat g,q)$ near $p$. If $\beta(p)=0$, then
we have $\gamma(p)>0$ and
\[
p\in\deltanbd{\what
S}{\delta_3}\cap\what\Omega\cap\Gamma\subset\what\Theta_1\cap\Gamma_0\sm\overline{\what\Omega_3}.
\]
Thus, near $p$, $\gamma=\gamma_0\circ\Upsilon$ is of class
$\strplshclass^\infty(\hat g,q)$ and $\beta$ is of class
$\plshclass^\infty(\hat g,q)$, and hence $\alpha$ is again of
class $\strplshclass^\infty(\hat g,q)$ near~$p$.

\textbf{Verification of (v).} Given $\nu\in\set{1,\dots,m}$, the
\fn $\beta-\hat\rho_\nu$ is of class $\strplshclass^\infty(\hat
g,q)$ on a \nbd of $\what D_\nu\cap\what\Theta_1\sm\deltanbd
Z{\delta_2}\supset\deltanbd{\what S}{\delta_3}\cap\what
D_\nu\cap\what\Omega\sm\deltanbd Z{\delta_2}$. Since $\gamma$ is
of class $\plshclass^\infty(\hat g,q)$ on
$\Gamma\cap\deltanbd{\what S}{\delta_3}$, the condition~(v) holds.

\textbf{Verification of (vi).} Given $\nu\in\set{1,\dots,m}$, we
have
\[
\what\Omega\cap\what D_\nu\sm\left[\deltanbd{\what
S}{\delta_2}\cup\deltanbd Z{\delta_1}\right]\subset\Gamma\cap\what
U\cap\what\Omega\cap\what D_\nu.
\]
Moreover, on $\Gamma\cap\what U\cap\what\Omega\cap\what D_\nu$, we
have
$\alpha-\hat\rho_\nu=(\beta+R_0\hat\rho_\nu)+(\gamma-(R_0+1)\cdot\hat\rho_\nu)$;
the sum of a \cinf \str \plsh \fn and a \cinf \strg $q$-\cvx
\fnns. Thus we get the property~(vi).

\textbf{Verification of (vii).} We have
\[
\what\Omega\sm\left[\deltanbd{\what S}{\delta_4}\cup\deltanbd
Z{\delta_1}\right]\subset\Gamma\cap\what\Omega\sm\deltanbd{\what
S}{\delta_4}\subset\Gamma\cap\what\Omega\cap\what\Omega_4.
\]
Therefore, on this set, we have $\alpha\geq\gamma>\theta$.

\textbf{Verification of (viii).} Recall that we extended the list
of \fns $\seq\rho\nu$ so that $\what\Omega\Subset\Theta_0\Subset
D_1\cup\cdots\cup D_m$. Hence the conditions (i),~(iv), and (vi)
give~(viii).
\end{pf}

\section{A uniformly quasi-\plsh exhaustion function}\label{Exh fn levi bdd sect}

To obtain the desired \exh \fnns, we will add to the \fn produced
in Proposition~\ref{q-cvx away from cpt not exhaustive prop} an
\exh \fn which is uniformly quasi-\plsh and which is locally
constant near the \cpt \ircomps of the lifting of the \cpt \anal
set. Such a \fn is produced in the following proposition which is
the main goal of this section:

\begin{prop}\label{Exh fn levi bdd below prop}
Let $(X,g)$ be a \con reduced Hermitian \cpx space, let $Y$ be a
\cpt \anal subset of~$X$, and, for each $\nu=1,\dots,m$, let
$D_\nu$ be a \rel \cpt open subset of $X$ and let $\rho_\nu$ be a
\cinf \str \plsh \fn on a \nbd of $\overline D_\nu$. Then there
exist constants $R>0$  and $\delta>0$ \stns, for every $\eta>0$,
for every set $E\subset\R$ with $|s-t|\geq 4\eta$ for all $s,t\in
E$ with $s\neq t$, for every set $F\subset\R$ with
$[0,\infty)\subset\deltanbdg F\eta{\R}$, for every \con covering
space $\Upsilon\colon\what X\to X$, and for every \anal set~$Z$
which is equal to a union of \ircomps of $\what Y=\Upsilon\inv(Y)$
and which has only \cpt \concompsns, there exists a \cinf positive
\exh \fn~$\tau$ on~$\what X$ with the following properties:
\begin{enumerate}
\item[(i)] The \fn $\tau+R\cdot(1+\eta)\cdot\rho_\nu\circ\Upsilon$
is \cinf \str \plsh on $\Upsilon\inv(D_\nu)$ for $\nu=1,\dots,m$;

\item[(ii)] For $\hat g=\Upsilon^*g$, $\tau$ is locally constant
on $\deltanbdg Z\delta{\hat g}$; and

\item[(iii)] We have $\tau(Z)\subset F\sm\deltanbdg E\eta{\R}$.

\end{enumerate}
\end{prop}

\begin{lem}\label{EF-lemma}
Let $\eta>0$, let $E\subset\R$ with $|s-t|\geq 4\eta$ for all
$s,t\in E$ with $s\neq t$, and let $F\subset\R$ with
$[0,\infty)\subset\deltanbdg F\eta{\R}$. Then, for each
$a\in[0,\infty)$, there exists a number $b\in[0,6\eta)$ \st
$a+b\in F\sm\deltanbdg E\eta{\R}$.
\end{lem}
\begin{pf}
For some $c\in F$, we have $|a+\eta-c|<\eta$ and hence
$a<c<a+2\eta$. If $c\notin\deltanbdg E\eta{\R}$, then we may set
$b=c-a$. If $|c-t|<\eta$ for some $t\in E$, then $t+2\eta>0$ and
hence, for some $d\in F$, we have $|t+2\eta-d|<\eta$. It follows
that $d\notin\deltanbdg E\eta{\R}$, since $\eta<|t-d|$ and, for
each $s\in E\sm\set{t}$,  $|s-d|\geq
|s-t|-|t-d|>4\eta-3\eta=\eta$. Moreover,
\[
a\leq c<t+\eta<d<t+3\eta<c+4\eta<a+6\eta.
\]
Thus $b=d-a$ has the required properties.
\end{pf}

\begin{lem}\label{Lipschitz fn on R with good growth}
Let $\seq s\nu$ and $\seq t\nu$ be sequences of nonnegative
numbers with $t_\nu\to\infty$ as $\nu\to\infty$. Then there exists
a Lipschitz \cont \fn $\chi$ on $\R$ \st
\begin{enumerate}
\item[(i)] We have $\chi>2$ and $0\leq\chi'\leq 1/2$;

\item[(ii)] We have $\chi(t)\to\infty$ as $t\to\infty$;

\item[(iii)] For all $s,t\in\R$ with $|s-t|<1$, we have
\(\chi(s)<2\cdot\chi(t)\); and

\item[(iv)] For all $\nu\in\N$ and all $t\in\R$ with
$t_\nu-1<t<s_\nu+t_\nu+1$, we have
\[
4\inv\cdot\chi(t)<\chi(t_\nu)<4\cdot\chi(t).
\]
\end{enumerate}
\end{lem}
\begin{pf}
Reordering if necessary, we may assume without loss of generality
that $\seq t\nu$ is nondecreasing. We may choose a subsequence
$\seq t{\nu_j}_{j=1}^\infty$ \st $t_{\nu_1}=t_1$ and \stns, for
each $j>1$, $t_{\nu_j}>s_\nu+t_\nu$ for all $\nu\leq\nu_{j-1}$ and
$t_{\nu_j}>t_{\nu_{j-1}}+2$. We now let $\chi$ be the piecewise
linear \cont \fn for which $\chi\equiv 3$ on $(-\infty,t_{\nu_1}]$
and, for $j>1$, $\chi(t_{\nu_j})=j+2$ and $\chi''\equiv 0$ on
$(t_{\nu_{j-1}},t_{\nu_j})$.

For each $\nu$, there is a unique index~$j$ with
$t_\nu\in[t_{\nu_{j-1}},t_{\nu_j})$ and hence
$j+1\leq\chi(t_\nu)<j+2$. If $s_\nu+t_\nu\leq t_{\nu_j}$, then
\[
\chi(s_\nu+t_\nu)\leq\chi(t_{\nu_j})=j+2\leq\frac{j+2}{j+1}\chi(t_\nu).
\]
Suppose $s_\nu+t_\nu>t_{\nu_j}$. Since $t_\nu<t_{\nu_j}$, we have
$\nu<\nu_j$ and, therefore, $t_{\nu_{j+1}}>s_\nu+t_\nu$ by the
choice of $t_{\nu_{j+1}}$. Thus
\[
\chi(s_\nu+t_\nu)\leq\chi(t_{\nu_{j+1}})=j+3\leq\frac{j+3}{j+1}\chi(t_\nu).
\]
Thus $\chi(s_\nu+t_\nu)\leq 2\chi(t_\nu)$ for each~$\nu$.

Since $\chi$ is piecewise linear, for each $j>1$ and each
$t\in(t_{\nu_{j-1}},t_{\nu_j})$, we have
\[
\chi'(t)=\frac{1}{t_{\nu_j}-t_{\nu_{j-1}}}<\frac 12.
\]
So (i) holds.

If $|s-t|<1$, then, since $\chi>2$ and $0\leq\chi'\leq 1/2$, we
have
\[
\chi(s)\leq\chi(t+1)\leq\chi(t)+\frac 12<2\cdot\chi(t).
\]
Hence~(iii) holds.

Finally, if $\nu\in\N$ and $t\in(t_\nu-1,s_\nu+t_\nu+1)$, then
$|s-t|<1$ for some $s\in[t_\nu,s_\nu+t_\nu]$. Therefore, by the
above,
\[
\frac 12\cdot\chi(t)<\chi(s)\leq\chi(s_\nu+t_\nu)\leq
2\cdot\chi(t_\nu)
\]
and
\[
\chi(t_\nu)\leq\chi(s)<2\cdot\chi(t).
\]
The property~(iv) now follows.
\end{pf}

\begin{pf*}{Proof of Proposition~\ref{Exh fn levi bdd below prop}}
We may fix a nonempty \rel \cpt \nbd $\Omega$ of $Y$ containing
$\overline D_\nu$ for $\nu=1,\dots,m$.  According to
Proposition~\ref{Metric modify to complete prop}, there exists a
Hermitian metric $g'\geq g$ on $X$ \st $g'=g$ on $\Omega$ and
$x\mapsto\dist_{g'}(x,p)$ is an \exh \fnns. If we find a
$\delta>0$ which works for this metric~$g'$ and which satisfies
$\delta<\dist_g(Y,X\sm\Omega)\leq\dist_{g'}(Y,X\sm\Omega)$, then
the $\delta$-\nbds in coverings of liftings of subsets of~$Y$ will
be the same for both metrics and hence this $\delta$ will also
work for the metric~$g$. Thus we may assume without loss of
generality that $x\mapsto\dist_g(x,p)$ is an \exh \fn for each
point $p\in X$.

We may first choose $\delta_0>0$ \st $\delta_0<1$,
$\delta_0<\dist_g(Y,X\sm\Omega)$, and, for each point
$p\in\overline\Omega$, $B_g(p,\delta_0)$ is contained in some
contractible open set. According to Lemma~\ref{Uniform dist
between comps in cover lem}, we may also choose $\delta_0$ so that
any two disjoint \ircomps of the lifting of $Y$ to a \con covering
space are of distance~$>\delta_0$ with respect to the lifting
of~$g$. We may now choose a number~$\delta$ with
$0<\delta<\delta_0/4$, a covering $V_1,\dots,V_k$ of
$\overline\Omega$ by finitely many \rel \cpt \con open subsets of
$X$ each of which is of diameter~$<\delta$ and meets $\Omega$, and
nonnegative \cinf \fns $\seq\alpha i_{i=1}^k$ \st
$\supp\alpha_i\subset V_i$ for each~$i$ and $\sum\alpha_i\equiv 1$
on a \nbd~$\Theta$ of $\overline\Omega$ in $X$.

Suppose now that $\eta>0$, $E\subset\R$ with $|s-t|\geq 4\eta$ for
all $s,t\in E$ with $s\neq t$, $F\subset\R$ with
$[0,\infty)\subset\deltanbdg F\eta{\R}$, $\Upsilon\colon\what X\to
X$ is a \con covering space,  $Z$ is an \anal subset of $X$ which
is equal to a union of \ircomps of $\what Y=\Upsilon\inv(Y)$ and
which has only \cpt \concompsns, and $\hat g=\Upsilon^*g$. Fixing
a point $O\in X$, Lemma~\ref{complete on covering lem} implies
that the \fn $r\colon x\mapsto\dist_{\hat g}(x,O)$ is an \exh
\fnns.  In order to produce the \fn~$\tau$, we may assume that the
covering is infinite since, for a finite covering, an \exh \fn
which is equal to a suitable constant on the lifting of $\Omega$
will work. Let $\set{Z_\sigma}_{\sigma\in S}$ be the distinct
\concomps of~$Z$. Observe that $\overline{\deltanbdg {\what
Y}{\delta_0}{\hat
g}}\subset\what\Omega\equiv\Upsilon\inv(\Omega)$.

By the choice of $\delta_0$ and $\delta$, each of the \cpt sets
$\overline V_i$ is contained in a contractible open set and is
therefore evenly covered by~$\Upsilon$. For each $i$, we let
$\set{V\ssp\nu_i}_{\nu=1}^\infty$ be the \concomps of $\what
V_i=\Upsilon\inv(V_i)$. The collection $\set{V\ssp\nu_i}$ is
locally finite in $\what X$ and (again, by the choice of
$\delta_0$ and $\delta$) each set $V\ssp\nu_i$ is of
diameter~$<\delta$ (with respect to the metric $\dist_{\hat
g}(\cdot,\cdot)$). For each $i=1,\dots,k$ and each
$\nu=1,2,3,\dots$, we define the nonnegative \cinf
\fn~$\alpha\ssp\nu_i$ with \cpt support in $V\ssp\nu_i$ by
\[
\alpha\ssp\nu_i(x)=\left\{
\begin{aligned}
\alpha_i(\Upsilon(x))&\qquad\text{if }x\in V\ssp\nu_i\\
0&\qquad\text{if }x\in \what X\sm V\ssp\nu_i
\end{aligned}
\right.
\]
We have $\sum_{i,\nu}\alpha\ssp\nu_i\equiv 1$ on the set
$\what\Theta=\Upsilon\inv(\Theta)$. By Lemma~\ref{Lipschitz fn on
R with good growth}, there exists a Lipschitz \cont \fn
$\chi\colon\R\to\R$ \st $\chi>2$, $0\leq\chi'\leq 1/2$,
$\chi(t)\to\infty$ as $t\to\infty$, \(\chi(s)<2\cdot\chi(t)\)
whenever $s,t\in\R$ with $|s-t|<1$, and
\[
\frac 14\cdot\chi(t)<\chi\bigl[\dist_{\hat
g}(O,Z_\sigma)\bigr]<4\cdot\chi(t)
\]
whenever $\sigma\in S$ and $\dist_{\hat
g}(O,Z_\sigma)-1<t<\diam_{\hat g}(Z_\sigma)+\dist_{\hat
g}(O,Z_\sigma)+1$. Finally, for each $\sigma\in S$, let $M_\sigma$
be the set of pairs of indices $(i,\nu)$ \st $1\leq i\leq k$,
$\nu\in\N$, and $V\ssp\nu_i\cap\deltanbdg{Z_\sigma}\delta{\hat
g}\neq\emptyset$. We also let~$M$ be the set of pairs of indices
$(i,\nu)$ \st $1\leq i\leq k$, $\nu\in\N$, and
$V\ssp\nu_i\cap\deltanbdg Z\delta{\hat g}=\emptyset$. Observe that
the sets $\seq M\sigma$ are mutually disjoint and disjoint
from~$M$. For, if $V\ssp\nu_i$ meets $\deltanbdg
{Z_\sigma}\delta{\hat g}$ and $\deltanbdg {Z_{\sigma'}}\delta{\hat
g}$, then, since $\diam_{\hat
g}\left(V\ssp\nu_i\right)<\delta<\delta_0/4$, we have $\dist_{\hat
g}(Z_\sigma,Z_{\sigma'})<\delta_0$. Hence $Z_\sigma\cap
Z_{\sigma'}\neq\emptyset$ and, therefore, $\sigma=\sigma'$.

We may now define a \cinf \fn~$\alpha$ on $\what X$ by
\[
\alpha=\sum_{\sigma\in S}\sum_{(i,\nu)\in
M_\sigma}\alpha\ssp\nu_i\cdot\chi\bigl[\dist_{\hat
g}(O,Z_\sigma)\bigr]+ \sum_{(i,\nu)\in
M}\alpha\ssp\nu_i\cdot\chi\left[\dist_{\hat
g}\left(O,V\ssp\nu_i\right)\right].
\]
We note that the previous remarks imply that the above sum is
locally finite. Moreover, each pair $(i,\nu)$ appears exactly once
and we have $\alpha>2$ on~$\what\Theta=\Upsilon\inv(\Theta)$. We
will show that the \fn $\beta=\log\alpha$ on $\what\Theta$ has the
following properties:

\begin{enumerate}

\item[(\ref{Exh fn levi bdd below prop}.1)] On the $\delta$-\nbd
$\deltanbdg Z\delta{\hat
g}\subset\what\Omega=\Upsilon\inv(\Omega)$ of~$Z$, $\beta$ is
locally constant;

\item[(\ref{Exh fn levi bdd below prop}.2)] The \fn $\beta$
exhausts $\overline{\what\Omega}$; and

\item[(\ref{Exh fn levi bdd below prop}.3)] There is a constant
$R_0>0$ which is independent of the choice of the covering space
$\what X$ and of the \anal subset $Z$ and for which the \fn
$\beta+R_0\cdot\rho_\nu\circ\Upsilon$ is \cinf \str \plsh on
$\Upsilon\inv(D_\nu)$ for $\nu=1,\dots,m$.

\end{enumerate}
For the property~(\ref{Exh fn levi bdd below prop}.1), we need
only observe that, for each $\sigma\in S$, we have
$\alpha\ssp\nu_i\equiv 0$ on~$\deltanbdg {Z_\sigma}\delta{\hat g}$
for each $(i,\nu)\not\in M_\sigma$. Thus, on $\deltanbd
{Z_\sigma}\delta$,
\[
\alpha=\sum_{(i,\nu)\in
M_\sigma}\alpha\ssp\nu_i\cdot\chi\biggl[\dist_{\hat
g}(O,Z_\sigma)\biggr]\equiv\chi\biggl[\dist_{\hat
g}(O,Z_\sigma)\biggr].
\]

For (\ref{Exh fn levi bdd below prop}.2), we first observe that,
for each index $\sigma\in S$, each pair~$(i,\nu)\in M_\sigma$, and
each point $x\in V\ssp\nu_i$, we have
\[
\dist_{\hat g}(O,Z_\sigma)-2\delta<r(x)=\dist_{\hat
g}(O,x)<\diam_{\hat g}(Z_\sigma)+\dist_{\hat
g}(O,Z_\sigma)+2\delta.
\]
For each pair~$(i,\nu)\in M$ and each point $x\in V\ssp\nu_i$, we
have
\[
\dist_{\hat g}(O,V\ssp\nu_i)<r(x)<\dist_{\hat
g}(O,V\ssp\nu_i)+\delta.
\]
Therefore, by the choice of $\chi$, we have
\(4\inv\cdot\chi(r)\leq\alpha\leq 4\cdot\chi(r)\) on
$\what\Theta$. In particular, $\alpha$ and $\beta$ exhaust
$\overline{\what\Omega}$, so the property~(\ref{Exh fn levi bdd
below prop}.2) holds.

Given a point $p\in\Theta$, we may choose a proper local \holo
model $(U,\Phi,U')$, a Hermitian metric~$g'$ on~$U'$, \cinf \fns
$\set{\alpha_i'}_{i=1}^k$ on $U'$, and a constant~$A>0$ \st $p\in
U\subset\Theta$, $\Phi^*g'=g$ on~$U$, and, for each $i=1,\dots,k$,
$\Phi^*\alpha_i'=\alpha_i$ on~$U$, and $|d\alpha_i'|_{g'}\leq A$
and $-A\cdot g'\leq\lev{\alpha_i'}\leq A\cdot g'$ on $U'$. We may
also assume that there is a constant $B>0$ \stns, whenever $U\cap
D_\nu\neq\emptyset$, there is a \cinf \str \plsh \fn $\rho_\nu'$
on $U'$ with $\rho_\nu=\rho_\nu'\circ\Phi$ on $U$ and
$B\cdot\lev{\rho_\nu'}\geq\left(16kA+256k^2A^2+1\right)\cdot g'$
on $U'$.

For each point $a\in\what U\equiv\Upsilon\inv(U)$, there is a
finite set $N$ of pairs of indices $(i,\nu)$ \st $a\in V\ssp\nu_i$
if and only if $(i,\nu)\in N$. In particular, for distinct
elements $(i,\nu),(j,\mu)\in N$, we have $i\neq j$; so $\#N\leq
k$. Thus, on some \rel \cpt \nbdns~$W$ of~$a$ in
$\bigcap_{(i,\nu)\in N}V\ssp\nu_i\cap\what U$, we have
\[
\alpha=\sum_{\sigma\in S}\sum_{(i,\nu)\in N\cap
M_\sigma}\alpha\ssp\nu_i\cdot\chi\bigl[\dist_{\hat
g}(O,Z_\sigma)\bigr]+ \sum_{(i,\nu)\in N\cap
M}\alpha\ssp\nu_i\cdot\chi\left[\dist_{\hat
g}\left(O,V\ssp\nu_i\right)\right].
\]
Setting $\Psi=\Phi\circ\Upsilon\restrict W$ and
\[
\alpha'=\sum_{\sigma\in S}\sum_{(i,\nu)\in N\cap
M_\sigma}\alpha_i'\cdot\chi\bigl[\dist_{\hat g}(O,Z_\sigma)\bigr]+
\sum_{(i,\nu)\in N\cap M}\alpha_i'\cdot\chi\left[\dist_{\hat
g}\left(O,V\ssp\nu_i\right)\right]\in\cinf(U'),
\]
we get $\Psi^*\alpha'=\alpha$ and $\Psi^*g'=\hat g$ on~$W$ and,
for each point $z=\Psi(x)$ with $x\in W$, we have
\begin{align*}
|(d\alpha')_z|_{g'}&=\left|\sum_{\sigma\in S}\sum_{(i,\nu)\in
N\cap M_\sigma}(d\alpha_i')_z\cdot\chi\bigl[\dist_{\hat
g}(O,Z_\sigma)\bigr]+ \sum_{(i,\nu)\in N\cap
M}(d\alpha_i')_z\cdot\chi\left[\dist_{\hat
g}\left(O,V\ssp\nu_i\right)\right]\right|_{g'}\\
&\leq\sum_{\sigma\in S}\sum_{(i,\nu)\in N\cap
M_\sigma}\left|(d\alpha_i')_z\right|_{g'}\cdot\chi\bigl[\dist_{\hat
g}(O,Z_\sigma)\bigr]\\
&\qquad\qquad+\sum_{(i,\nu)\in N\cap
M}\left|(d\alpha_i')_z\right|_{g'}\cdot\chi\left[\dist_{\hat
g}\left(O,V\ssp\nu_i\right)\right]\\
&\leq\sum_{\sigma\in S}\sum_{(i,\nu)\in N\cap
M_\sigma}\left|(d\alpha_i')_z\right|_{g'}\cdot
4\cdot\chi(r(x))+\sum_{(i,\nu)\in N\cap
M}\left|(d\alpha_i')_z\right|_{g'}\cdot 4\cdot\chi(r(x))\\
&\leq 4kA\cdot \chi(r(x)).
\end{align*}
Similarly, for each point $z=\Psi(x)$ with $x\in W$ and each
tangent vector $v\in T^{(1,0)}_zU'$, we have
\[
\left|\lev{\alpha'}(v,v)\right|\leq
4kA\cdot|v|^2_{g'}\cdot\chi(r(x)).
\]
Setting $\beta'=\log\alpha'$, we get $\beta=\beta'\circ\Psi$ on
$W$ and, for each point $z=\Psi(x)$ with $x\in W$ and each tangent
vector $v\in T^{(1,0)}_zU'$, we get
\begin{align*}
\lev{\beta'}(v,v)&=\frac 1{\alpha(x)}\cdot\lev{\alpha'}(v,v)-\frac
1{(\alpha(x))^2}|\partial\alpha'(v)|^2\\
&\geq-\frac {4kA\cdot\chi(r(x))}{\alpha(x)}\cdot|v|^2_{g'}-\frac
{16k^2A^2\cdot(\chi(r(x)))^2}{(\alpha(x))^2}\cdot|v|^2_{g'}\\
&\geq-\left(16kA+256k^2A^2\right)\cdot|v|^2_{g'}.
\end{align*}
Hence, whenever $U\cap D_\nu\neq\emptyset$, the \fn
$\beta'+B\rho_\nu'$ will be \str \plsh on a \nbd of $\Psi(W)$ in
$U'$ and will satisfy
$\beta+B\rho_\nu\circ\Upsilon=(\beta'+B\rho_\nu')\circ\Psi$ on
$W$. Covering $\overline\Omega$ by finitely many such \nbdsns~$U$,
we see that we may choose a constant $R_0>0$ so that the
property~(\ref{Exh fn levi bdd below prop}.3) holds for all
covering spaces $\what X$ and all choices of $Z$.

We now modify $\beta$ to get a \fn with values in $F\sm\deltanbdg
E\eta{\R}$ on $Z$. According to Lemma~\ref{EF-lemma}, for each
$\sigma\in S$, we may choose a number $b_\sigma\in[0,6\eta)$ \st
\[
c_\sigma=\log\left(\chi\left[\dist_{\hat
g}(O,Z_\sigma)\right]\right)+b_\sigma\in F\sm\deltanbdg E\eta{\R}.
\]
We may then define a \cinf \fnns~$\gamma$ on $\what X$ by the
locally finite sum
\[
\gamma=\sum_{\sigma\in S}\sum_{(i,\nu)\in
M_\sigma}\alpha\ssp\nu_i\cdot b_\sigma.
\]
Easier versions of the previous arguments show that~$\gamma\equiv
b_\sigma$ on $\deltanbdg Z\delta{\hat g}\subset\what\Omega$ for
each $\sigma\in S$, and, for some constant $R_1>0$ independent of
the choice of the covering $\what X$, the \anal set $Z$, and the
number $\eta$, the \fn $\gamma+R_1\eta\cdot\rho_\nu\circ\Upsilon$
will be \cinf \str \plsh on $\Upsilon\inv(D_\nu)$ for
$\nu=1,\dots,m$. Setting $R=\max(R_0,R_1)$, it now follows that
any positive \cinf \exh \fn $\tau$ on $X$ which is equal to
$\beta+\gamma$ on $\what\Omega$ has the required properties.
\end{pf*}

\section{The main result}\label{Main result section}

Theorem~\ref{main theorem from intro} is an immediate consequence
of the following theorem which is proved in this section:
\begin{thm}\label{Main thm generl version}
Let $(X,g)$ be a \con reduced Hermitian \cpx space; let $q$ be a
positive integer; let $Y$ be a \cpt \anal subset of
dimension~$\leq q$; let~$C$~and~$S$ be \anal subsets of~$Y$ \st
$C$ is a union of \ircomps of $Y$, $\dim S<q$, $S$ contains $\sing
Y$ as well as every \ircomp of $Y$ of dimension~$<q$, and $Y\sm
(C\cup S)$ is Stein; let $\seq D\nu_{\nu=1}^m$ be \rel \cpt open
subsets of~$X$; let $\rho_\nu$ be a \cinf \str \plsh \fn on a \nbd
of $\overline D_\nu$ for each $\nu=1,\dots,m$; and let
$\epsilon>0$. Then, for every choice of constants
$\delta_1$,~$\delta_2$,~$\delta_3$,~and $\delta_4$ with
$\epsilon\gg\delta_4\gg\delta_3\gg\delta_2\gg\delta_1>0$ (i.e.~one
must choose $\delta_4$ sufficiently small relative to~$\epsilon$,
$\delta_3$ sufficiently small relative to $\delta_4$, and so on),
there exists a \nbd $\Omega$ of $Y$ in $X$ \stns, for every
$\eta>0$, for every set $E\subset\R$ with $|s-t|\geq 4\eta$ for
all $s,t\in E$ with $s\neq t$, for every set $F\subset\R$ with
$[0,\infty)\subset\deltanbdg F\eta{\R}$, for every \con covering
space $\Upsilon\colon\what X\to X$, for every union~$Z$ of
\ircomps of $\what Y=\Upsilon\inv(Y)$ which has only \cpt
\concomps and which contains $\what C=\Upsilon\inv(C)$ as well as
every $q$-dimensional \cpt \ircomp of $\what Y$ (i.e. $Z$ contains
every \cpt \ircomp of $\what Y$ not contained in $\what
S=\Upsilon\inv(S)$), and for every positive \cont \fn $\theta$ on
$\what X$, there exists a positive \cinf \exh \fn $\vphi$ on
$\what X$ with the following properties relative to the Hermitian
metric $\hat g=\Upsilon^*g$, the sets
$\what\Omega=\Upsilon\inv(\Omega)$ and $\what
D_\nu=\Upsilon\inv(D_\nu)$ for $\nu=1,\dots,m$, and the \fns
$\hat\rho_\nu=\rho_\nu\circ\Upsilon\colon\what D_\nu\to\R$ for
$\nu=1,\dots,m$:
\begin{enumerate}

\item[(i)] On $\deltanbd Z{\delta_1}\cap\what\Omega$, $\vphi$ is
locally constant with
\[
G=\vphi\left(\deltanbd
Z{\delta_1}\cap\what\Omega\right)=\vphi(Z)\subset F\sm\deltanbdg
E\eta{\R};
\]

\item[(ii)] On $\deltanbd{\what S}{\delta_3}\cap\what\Omega$,
$\vphi$ is of class $\plshclass^\infty(\hat g,q)$;

\item[(iii)] On the set $\setof{p\in\deltanbd {\what
S}{\delta_3}\cap\what\Omega}{\vphi(p)\notin G}$, $\vphi$ is of
class $\strplshclass^\infty(\hat g,q)$;

\item[(iv)] For each $\nu=1,\dots,m$, the \fn $\vphi-\hat\rho_\nu$
is of class $\strplshclass^\infty(\hat g,q)$ on a \nbd of
$\deltanbd{\what S}{\delta_3}\cap\what
D_\nu\cap\what\Omega\sm\deltanbd Z{\delta_2}$;

\item[(v)] For each $\nu=1,\dots,m$, the \fn $\vphi-\hat\rho_\nu$
is \cinf \strg $q$-\cvx on a \nbd of $\what
D_\nu\cap\what\Omega\sm\left[\deltanbd {\what
S}{\delta_2}\cup\deltanbd Z{\delta_1}\right]$;

\item[(vi)] On $\what\Omega\sm\left[\deltanbd{\what
S}{\delta_4}\cup\deltanbd Z{\delta_1}\right]$, $\vphi>\theta$;

\item[(vii)] On $\setof{p\in\what\Omega}{\vphi(p)\notin G}$,
$\vphi$ is \cinf \strg $q$-\cvxns; and

\item[(viii)] On $\what\Omega$, $\vphi$ is of class
$\plshclass^\infty(q)$.

\end{enumerate}
\end{thm}
\begin{pf}
We may assume without loss of generality that $\deltanbd
Y\epsilon\Subset D_1\cup\cdots\cup D_m$. Applying
Proposition~\ref{Exh fn levi bdd below prop}, we may also assume
that, for some constant $R_0>1$, for every $\eta>0$, for every set
$E\subset\R$ with $|s-t|\geq 4\eta$ for all $s,t\in E$ with $s\neq
t$, for every set $F\subset\R$ with $[0,\infty)\subset\deltanbdg
F\eta{\R}$, for every \con covering space $\Upsilon\colon\what
X\to X$, and for every \anal set~$Z$ which is equal to a union of
\ircomps of $\what Y=\Upsilon\inv(Y)$ and which has only \cpt
\concompsns, there exists a \cinf positive \exh \fn~$\tau$
on~$\what X$ with the following properties:
\begin{enumerate}
\item[(\ref{Main thm generl version}.1)] The \fn
$\tau+R_0\cdot(1+\eta)\cdot\rho_\nu\circ\Upsilon$ is \cinf \str
\plsh on $\Upsilon\inv(D_\nu)$ for $\nu=1,\dots,m$; and

\item[(\ref{Main thm generl version}.2)] For $\hat g=\Upsilon^*g$,
$\tau$ is locally constant on $\deltanbdg Z\epsilon{\hat g}$ and
\(\tau\left(\deltanbdg Z\epsilon{\hat g}\right)\subset
F\sm\deltanbdg E\eta{\R}\).

\end{enumerate}
Note that, in the above, $\tau(\deltanbd Z\epsilon)=\tau(Z)$ since
$\tau$ is constant on the \con \nbd $\deltanbd {Z'}\epsilon$ of
each \concomp $Z'$ of $Z$ and $\deltanbd Z{\epsilon}$ is equal to
the union of all such \nbdsns.

According to Proposition~\ref{q-cvx away from cpt not exhaustive
prop}, if we choose constants
$\delta_1$,~$\delta_2$,~$\delta_3$,~and $\delta_4$ with
$\epsilon\gg\delta_4\gg\delta_3\gg\delta_2\gg\delta_1>0$, then
there exists a \rel \cpt \nbd $\Omega$ of $Y$ in $\deltanbd
Y{\delta_1}$ \stns, for every \con covering space
$\Upsilon\colon\what X\to X$, for every \anal set $Z$ which is
equal to a union of \ircomps of $\what Y\equiv\Upsilon\inv(Y)$ and
which contains $\what C=\Upsilon\inv(C)$ as well as every
$q$-dimensional \cpt \ircomp of $\what Y$, and for every positive
\cont \fn $\theta$ on $\what X$, there exists a nonnegative \cinf
\fn $\alpha$ on $\what X$ with the following properties relative
to the Hermitian metric $\hat g=\Upsilon^*g$, the sets
$\what\Omega=\Upsilon\inv(\Omega)$, $\what S=\Upsilon\inv(S)$, and
$\what D_\nu=\Upsilon(D_\nu)$ for $\nu=1,\dots,m$, and the \fns
$\hat\rho_\nu=\rho_\nu\circ\Upsilon\colon\what D_\nu\to\R$ for
$\nu=1,\dots,m$:
\begin{enumerate}

\item[(\ref{Main thm generl version}.3)] On $\deltanbd
Z{\delta_1}\cap\what\Omega$, $\alpha\equiv 0$;

\item[(\ref{Main thm generl version}.4)] On
$\left(\what\Omega\sm\deltanbd
Z{\delta_2}\right)\cup\left(\what\Omega\sm\left[\deltanbd{\what
S}{\delta_2}\cup\deltanbd Z{\delta_1}\right]\right)$, $\alpha>0$;

\item[(\ref{Main thm generl version}.5)] On $\deltanbd{\what
S}{\delta_3}\cap\what\Omega$, $\alpha$ is of class
$\plshclass^\infty(\hat g,q)$;

\item[(\ref{Main thm generl version}.6)] On the set
$\setof{p\in\deltanbd {\what
S}{\delta_3}\cap\what\Omega}{\alpha(p)>0}$, $\alpha$ is of class
$\strplshclass^\infty(\hat g,q)$;

\item[(\ref{Main thm generl version}.7)] For each $\nu=1,\dots,m$,
the \fn $\alpha-\hat\rho_\nu$ is of class
$\strplshclass^\infty(\hat g,q)$ on a \nbd of $\deltanbd{\what
S}{\delta_3}\cap\what D_\nu\cap\what\Omega\sm\deltanbd
Z{\delta_2}$;

\item[(\ref{Main thm generl version}.8)] For each $\nu=1,\dots,m$,
the \fn $\alpha-\hat\rho_\nu$ is \cinf \strg $q$-\cvx on a \nbd of
$\what D_\nu\cap\what\Omega\sm\left[\deltanbd{\what
S}{\delta_2}\cup\deltanbd Z{\delta_1}\right]$;

\item[(\ref{Main thm generl version}.9)] On
$\what\Omega\sm\left[\deltanbd{\what S}{\delta_4}\cup\deltanbd
Z{\delta_1}\right]$, $\alpha>\theta$; and

\item[(\ref{Main thm generl version}.10)] On
$\setof{p\in\what\Omega}{\alpha(p)>0}$, $\alpha$ is \cinf \strg
$q$-\cvxns.

\end{enumerate}

Suppose now that $\eta>0$, $E\subset\R$ with $|s-t|\geq 4\eta$ for
all $s,t\in E$ with $s\neq t$, $F\subset\R$ with
$[0,\infty)\subset\deltanbdg F\eta{\R}$, $\Upsilon\colon\what X\to
X$ is a \con covering space, $Z$ is a union of \ircomps of $\what
Y=\Upsilon\inv(Y)$ which has only \cpt \concomps and which
contains $\what C=\Upsilon\inv(C)$ as well as every
$q$-dimensional \cpt \ircomp of $\what Y$, and $\theta$ is a
positive \cont \fn on $\what X$. Setting $\hat g=\Upsilon^*g$,
$\what\Omega=\Upsilon\inv(\Omega)$, $\what S=\Upsilon\inv(S)$, and
$\what D_\nu=\Upsilon\inv(D_\nu)$ and
$\hat\rho_\nu=\rho_\nu\circ\Upsilon\colon\what D_\nu\to\R$ for
$\nu=1,\dots,m$, we get on $\what X$ a positive \cinf \exh
\fnns~$\tau$ and a nonnegative \cinf \fn $\alpha$ satisfying
(\ref{Main thm generl version}.1)--(\ref{Main thm generl
version}.10). It is then easy to check that the positive \cinf
\exh \fn
\[
\vphi\equiv\tau+R_0\cdot(2+\eta)\cdot\alpha
\]
has the required properties (i)--(viii).
\end{pf}

\bibliographystyle{amsalpha}

\end{document}